\documentclass[a4paper,parskip,11pt]{scrartcl}
\usepackage[utf8]{inputenc}
\usepackage[english]{babel} 
\usepackage{latexsym} 
\usepackage{pgfplots}
\usepackage{amsmath}
\usepackage{amsfonts}
\usepackage{amsthm}
\usepackage{amssymb}
\usepackage[top=1in,bottom=1in,right=1in,left=1in]{geometry}
\usepackage{graphicx}
\usepackage{nameref}
\usepackage{stmaryrd}
\usepackage{bbold} 
\usepackage{pgfplots}
\usepackage[normalem]{ulem}
\usepackage{abstract}
\usepackage{fancybox}
\usepackage{cases}

\usepackage{subfig}
\usepackage{pgfplots}
\pgfplotsset{compat=newest}
\usepackage{float}
\usepackage{caption}

\usepackage{mdwlist}
\usepackage{eurosym}
\usepackage{lastpage}
\usepackage{framed}
\usepackage{xcolor}
\usepackage[footsepline,plainfootsepline]{scrlayer-scrpage}
\usepackage{tcolorbox}
\usepackage{mathtools}
\usepackage{graphicx}
\usepackage[colorlinks=true,bookmarks=true,citecolor=black,linkcolor=black,pdftitle={Weak and strong interaction of excitations kinks in scalar parabolic equations}]{hyperref}
\usepackage{cleveref}
\usepackage{paralist}
\usepackage{multicol}

\usepackage{todonotes}

\definecolor{matlabblue}{rgb}{0.3010, 0.7450, 0.9330}
\definecolor{matlaborange}{rgb}{0.9100, 0.4100, 0.1700}
\definecolor{matlabyellow}{rgb}{0.9290, 0.6940, 0.1250}
\definecolor{matlabpurple}{rgb}{0.4940, 0.1840, 0.5560}

\def\Z{\mathbb{Z}} 
\def\N{\mathbb{N}} 
\def\R{\mathbb{R}} 
\def\calM{\mathcal{M}} 
\def\calK{\mathcal{K}} 
\def\calU{\mathcal{U}} 

\def\eps{\varepsilon}

\def\P{w}

\def\Rrm{\textup{\textrm{R}}}
\def\Lrm{\textup{\textrm{L}}}

\def\td{\tilde{\delta}}
\def\muf{f_0}

\renewcommand{\hat}{\widehat}

\def\calO{\mathcal{O}}

\def\ds{\displaystyle}
\def\rmd{\mathrm{d}}
\def\dr{\partial}
\def\vp{\varphi}
\def\ML{\mathcal{L}}
\def\MX{\mathcal{X}}
\def\MN{\mathcal{N}}

\def\lp{\left(}
\def\rp{\right)}

\def\lb{\left|}
\def\rb{\right|}
\def\lV{\left\Vert}
\def\rV{\right\Vert}

\newtheorem{theorem}{Theorem}
\newtheorem{lemma}[theorem]{Lemma}

\newtheorem{proposition}[theorem]{Proposition}
\newtheorem{remark}[theorem]{Remark}

\pdfoutput=1

\title{Weak and strong interaction of excitation kinks in scalar parabolic equations}
\author{A. Pauthier, J.D.M. Rademacher, D. Ulbrich}
\date{\today}

\begin{document}

\maketitle

\begin{abstract}
Motivated by studies of the Greenberg-Hastings cellular automata (GHCA) as a caricature of excitable systems, in this paper we study kink-antikink dynamics in the perhaps simplest PDE model of excitable media given by the scalar reaction diffusion-type $\theta$-equations for excitable angular phase dynamics. On the one hand, we qualitatively study geometric kink positions using the comparison principle and the theory of terraces. This yields the minimal initial distance as a global lower bound, a well-defined sequence of collision data for kinks- and antikinks, and implies that periodic pure kink sequences are asymptotically equidistant. On the other hand, we study metastable dynamics of finitely many kinks using weak interaction theory for certain analytic kink positions, which admits a rigorous reduction to ODE. By blow-up type singular rescaling we show that distances become ordered in finite time, and eventually diverge.
We conclude that diffusion implies a loss of information on kink distances so that the entropic complexity based on positions and collisions in the GHCA does not simply carry over to the PDE model. 
\end{abstract}


\section{Introduction}
Many spatially extended physical, chemical and biological systems form so-called excitable media, in which a 
supercritical perturbation from a stable equilibrium triggers an excitation that is transferred to its neighbours, followed 
by refractory return to the rest state. 

\begin{figure}[t]
\centering
\begin{tabular}{cc}
    \hspace{-0.5cm}\includegraphics[scale=0.43]{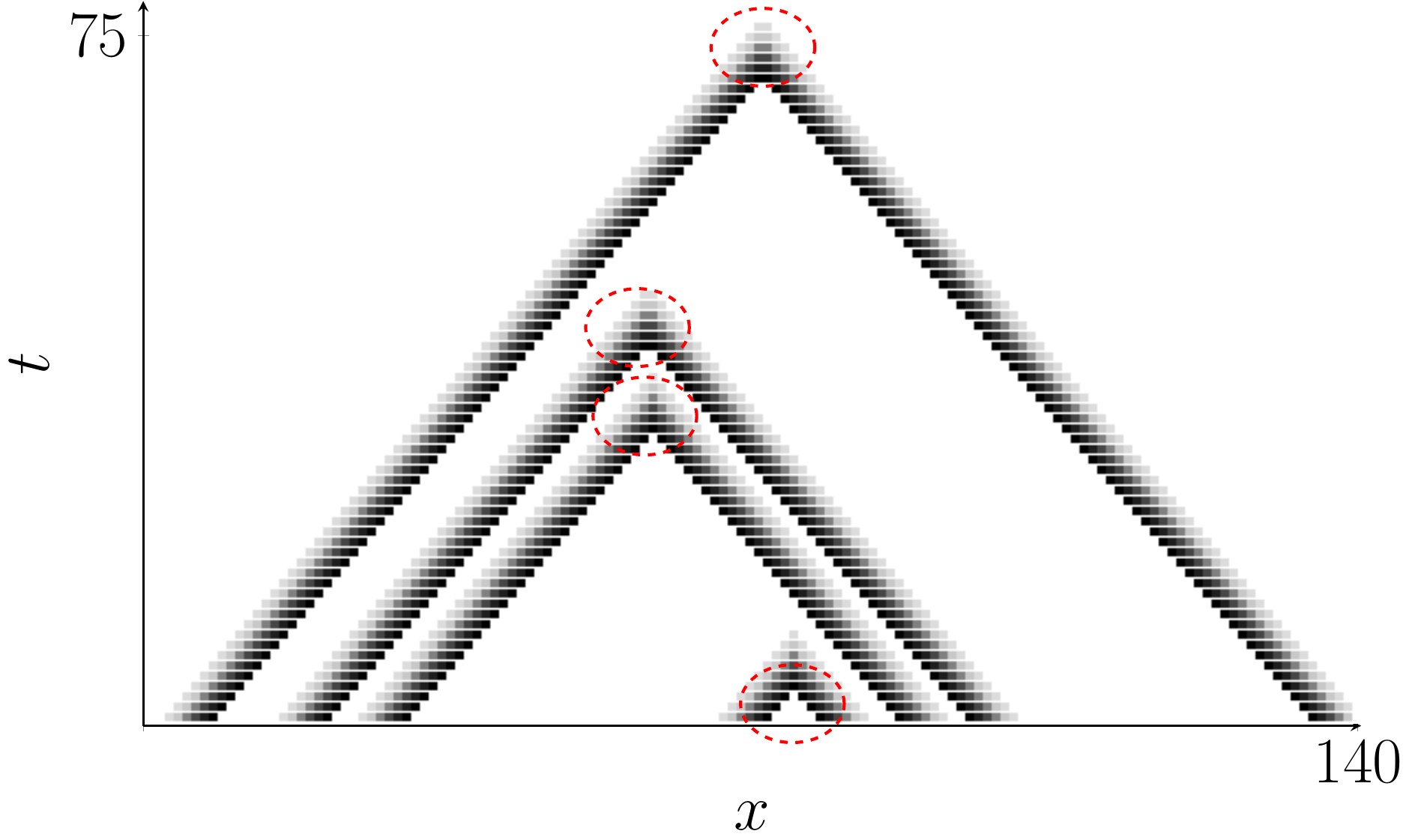}&
    \hspace{-0.6cm}\includegraphics[scale=0.43]{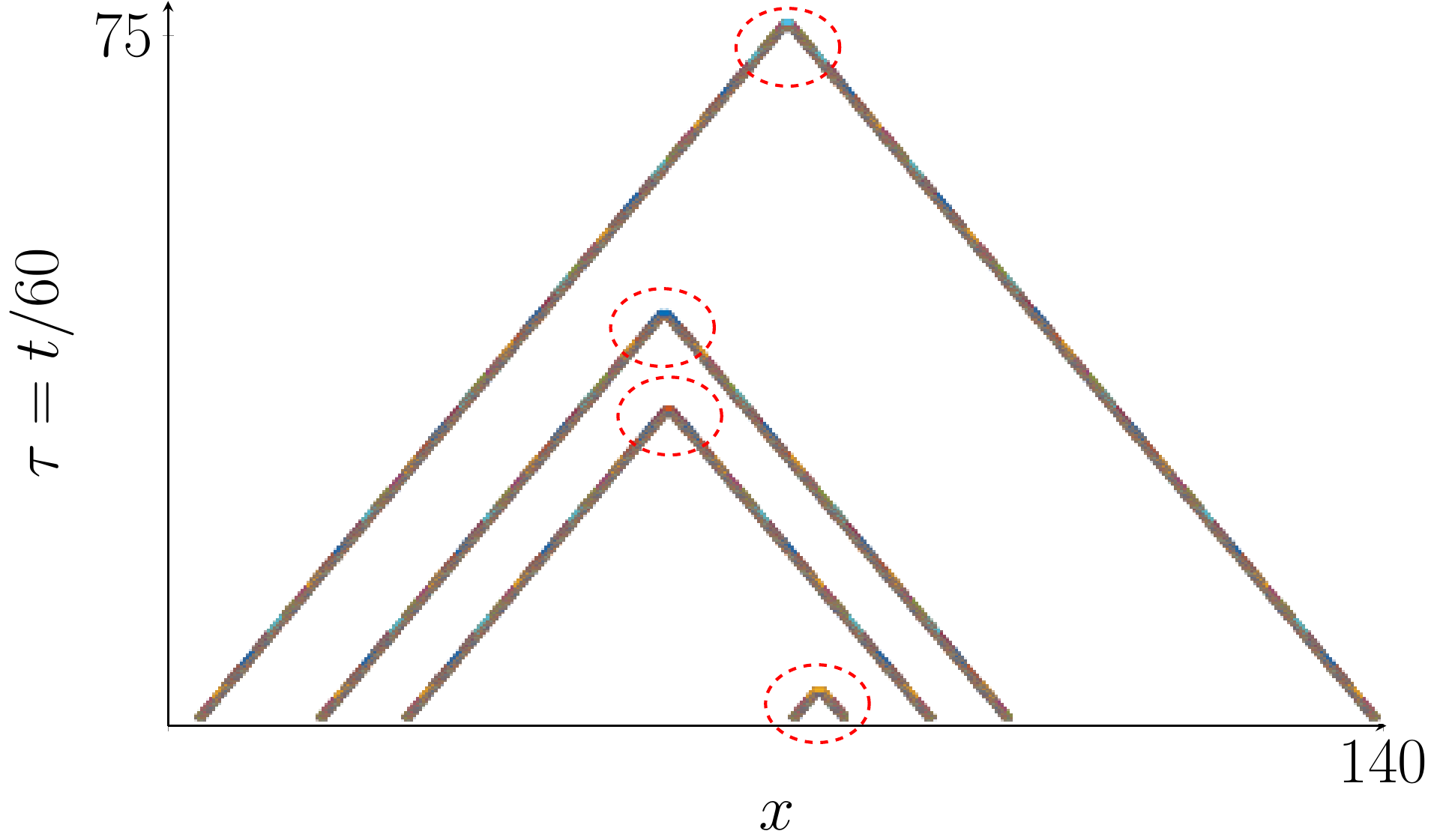}\\
    (a) & (b)
    \end{tabular}
\caption{Space-time plots of pulse positions for GHCA (a) and positions \eqref{e:pset} of \eqref{theta-eq} (b) with $\muf=0.05$ for initial data with four pairs of pulses and kinks/antikinks, respectively; marked are the associated annihilation events. In the right panel, time is rescaled such that a single front has unit speed as in the left panel. Some details on the numerics can be found in Appendix \ref{app:sim}.
}
\label{fig:GHCATheta}
\end{figure}

Perhaps the simplest dynamical system that realises a caricature excitable medium is the 1D Greenberg-Hastings 
cellular automata (GHCA) \cite{GreHas78,greenberg1978pattern}. Its key features are that spatial excitation loops embedded between rest 
states form local pulses that travel in one direction, and a counter-propagating pair of such local pulses annihilates, leading to a pure rest state locally, cf. Fig.~\ref{fig:GHCATheta} (a). 
Pulse trains can be generated from local pulses by placing these at arbitrary positions; they will maintain their relative distance 
until a possible annihilation. In fact, the topological entropy of the GHCA results from a Devaney-chaotic closed invariant subset of the non-wandering 
set that consists of colliding and annihilating local pulses \cite{DurSte91,kessebhmer2019dynamics}. From a dynamical 
systems viewpoint the non-wandering set is of prime relevance, and it turns out that for the GHCA it can be decomposed into 
invariant sets with different wave dynamics \cite{kessebhmer2019dynamics}. 

However, the modelling of excitable media is predominantly done by parabolic partial differential equations 
(PDE) and systems thereof. A paradigm is the famous FitzHugh-Nagumo (FHN) equation, derived from the Hodgin-Huxley 
model for nerve axons \cite{fitzhugh1961impulses,nagumo1962active,hodgkin1952quantitative}. 
A priori, for any such PDE model, the identification and description of pulse dynamics akin to GHCA is a formidable task and far from completely understood even in FHN. 
The PDE dynamics is in general also richer and parameter dependent, for instance self-replication of pulses has been numerically observed \cite{carter2020pulse,hayase2000self}, and rebound of pulses upon collision \cite{nishiura2000self}. 

While the existence of arbitrary local pulse trains and their annihilation is trivial in the GHCA, already the existence proof of a single pulse 
in FHN is not. It is well known that FHN possesses a homoclinic travelling wave solution that is 
spatially asymptotic to a stable rest state and thus corresponds to a single local pulse in GHCA; by spatial reflection the 
direction of motion can be reversed. However, there is no meaningful notion of a \emph{local} pulse since the spatial 
coupling by diffusion is effectively non-local through its infinite propagation speed. An analogue of pulse trains on the level of initial data is a 
superposition of multiple  travelling wave pulses placed at a distance from each other. In the dynamics of the PDE, diffusivity immediately 
couples these pulses, albeit in a weak form. In the past decades, a number of results have been obtained that rigorously relate positions of 
these `pulses' to an ordinary differential equation system (ODE) \cite{Ei2002,ZelikMielke}. However, to our knowledge there are no analytical results for collisions of excitation pulses in higher order equations or systems of PDE -- 
the closest result is the existence of an attracting invariant manifold for a sufficiently distant pair of counter-propagating pulses 
\cite{scheel2008colliding,wright2009separating}. The situation simplifies for scalar parabolic equations; in particular, the comparison principle for scalar parabolic equations allows to study collisions. We note that also energy methods can be used \cite{West2020}, even for the fourth order (non-excitable) Cahn-Hilliard equation~\cite{Scholtes2018}.

\medskip
In this paper we consider suitable scalar parabolic PDE for periodic phase dynamics as models for excitable media, and study similarities and 
difference between this and the local pulse dynamics of the GHCA. These so-called $\theta$-equations for oscillator phase dynamics 
\cite{ErmentroutRinzel,ScheelBellay} are given by  
\begin{equation}\label{theta-eq}
    \theta_t=\theta_{xx}+f(\theta),\qquad \theta(t,x)\in S^1 = \mathbb{R}/2\pi\mathbb{Z},\, x\in\mathbb{R}\,, t>0.
\end{equation}
For simplicity, we specify the nonlinearity as $f(\theta):=\cos(\theta+\Theta_0)+\muf$, where $\muf\in(0,1)$ and $\Theta_0\in(0,\pi)$ is chosen uniquely so that $f(0)=0$, although all results are valid in general for the excitable case of \cite{ErmentroutRinzel}. 

Equation \eqref{theta-eq} possesses the stable spatially homogeneous state $\theta\equiv 0$ and a right-moving (as well as left-moving upon reflection) travelling wave solution $\theta(t,x) = \check\varphi(\xi)$, $\xi=\pm x-ct$ with speed $c>0$.
The profile $\check\varphi(\xi)$ is asymptotic to $0\in S^1$ for $|\xi|\to \infty$, but with non-trivial winding number. 
This means that upon lifting $S^1$ to $\mathbb{R}$, the rest state 
maps to the sequence $2k\pi$, $k\in \mathbb{Z}$, and $\check\varphi$ maps to heteroclinic \emph{front} solutions 
$\varphi+2\pi k$ with $\varphi(\xi)\to 0$ as $\xi\to\infty$, and $\varphi(\xi)\to 2\pi$ as $\xi\to-\infty$.

The choice of $f$ explicitly relates \eqref{theta-eq} to the overdamped limit of the Sine-Gordon equation, which also arises as a phase 
field model for certain complex Ginzburg-Landau equations with broken gauge symmetry \cite[eqn.\ (91)]{AransonKramer}.
Following the terminology of the Sine-Gordon equation, we refer to these fronts as \emph{kinks}, and their left-moving spatial reflection as \emph{antikinks}. 

In the present context we view these as corresponding to a single local pulse in GHCA. We are thus concerned with the evolution of initial data built from kink- and antikink sequences, such as plotted in 
Fig.~\ref{fig:GHCATheta}. For any given solution $u(x,t)$ viewed in the lift to $\mathbb{R}$, we geometrically define the set of positions of potential kinks and antikinks through the phase $\pi$ intercepts as 
\begin{equation}\label{e:pset}
P(u(t)) := \{\xi\in \mathbb{R} |  \exists k\in\mathbb{Z}: u(\xi,t) = (2k+1)\pi \},
\end{equation}
which readily turns out to be a discrete set for $t>0$. We will discuss conditions under which this set meaningfully encodes kink or antikink 
positions also for large $t\gg 1$. A crucial ingredient is the theory of terraces.

Roughly speaking, a terrace is a superposition of finitely many fronts, each of them being a heteroclinic connection 
between two equilibria, such that asymptotically these fronts are separated, i.e., the distances between any two fronts 
diverge. To our knowledge, the terminology was introduced in \cite{DucrotGilettiMatano}, even though the notion 
was already present in the seminal paper \cite{Fife1977} for fronts with distinct speeds. Since then, it appeared that this 
notion is fundamental in the understanding of the long time behaviour of solutions of reaction-diffusion equations.
Pol\'{a}\v{c}ik proved in \cite{polacik2016propagating} under mild assumptions, in the context of scalar homogeneous semi-linear parabolic equations, that any front-like initial data eventually leads to a terrace profile. 
These results were extended to localized initial data in \cite{MatanoPolacik2}, leading to a pair of terrace profiles going in opposite directions. 
Risler independently proved similar results in the more general context of gradient flows in \cite{Risler} for systems. 

When all the considered fronts travel at asymptotically distinct speeds, the behaviour of the terrace profile, the 
convergence toward it (starting from front-like initial data, for example), as well as its stability are quite well understood, see  \cite{polacik2016propagating,LinSchecter} and references therein. On the other hand, if two consecutive fronts have the same asymptotic speed, the question is much more
intricate since the eventual separation is powered by weak interactions. In our context, kinks are right moving fronts, while antikinks are left moving fronts, all at the same speed. 

Hence, the model \eqref{theta-eq} is suited for a comparison with GHCA as it allows to study the combination of strong interactions, when a kink and an antikink collide, \emph{and} weak interactions when considering consecutive kinks. However, this observation already suggests a multi-scale nature, which turns out to imply a fundamental difference for the long time dynamics of excitation pulse positions.

We informally summarize our main results. Suited to unbounded kink-antikink initial data, we prove global well-posedness for unbounded data,  
and infer that the set of geometric positions \eqref{e:pset} consists of isolated points that lie on smooth curves except possibly at collisions.  
On the one hand, as a qualitative analysis, we rely on the comparison principle and roughly track positions which allows us to show that the minimal initial distance is a global lower bound for distances, and to abstractly identify collision times and locations. Moreover, as a model for the dynamics far from collisions, we show that under periodic boundary conditions the distances of neighboring kinks asymptotically equalise, and we numerically corroborate that this loss of information happens more broadly.
On the other hand, as a quantitative analysis, we derive the ODE in the weak interaction regime following \cite{Ei2002, Rossides2015} for finitely many kinks and analyse these in some detail, which extends the aforementioned qualitative terrace results. Our analysis relies on blow-up type singular rescaling and identifies the dynamics as being slaved to dynamics on a sphere at infinity, which, for instance, shows that distances become ordered in finite time. 
The latter is again hinting a loss of information through the dynamics: the memory of the initial distances is `washed out'. 
This suggests that the chaotic dynamics of GHCA and the entropic complexity based on positional dynamics is reduced by the weak interaction in the PDE, and we do not expect a topological entropy (in a suitable sense for the PDE) -- based on positions alone -- which resembles that of GHCA. 

This paper is organized as follows. In \S\ref{s:positions} we discuss the notions of positions and the initial data. Section~\ref{boundedmonotone} is devoted to qualitative and quantitative aspects for bounded monotone data which are composed solely of kinks (or antikinks). 
In \S\ref{boundedkinkantikink}, we focus on the annihilation process for initial data composed of  kinks and antikinks.
Unbounded initial data, i.e., infinitely many kinks and antikinks, as well as the long-time complexity are considered in \S\ref{sec:unbounded}.
Finally, in \S\ref{sec:disc} we conclude with a discussion. In the appendices we collect some technical proofs as well as notes on the numerical implementations.

\section{Kinks and anti-kinks and their positions}\label{s:positions}

In this section we introduce initial data which has imprinted positions of kinks and anti-kinks such that the setwise definition of `positions' $P(u(t))$ in \eqref{e:pset}, which is well-defined as long as $u(t)$ is defined, turns into meaningful individual positions for all time. 
In order to discuss this further, we first turn to the notion of kinks and anti-kinks in more detail.

The aforementioned regime $\muf\in (0,1)$ is termed excitable since the ODE for spatially constant data possesses a stable rest state and an unstable state which acts as a threshold for undergoing an `excitation loop', i.e., winding once through $S^1$. For $\muf>1$ the two equilibria have undergone a saddled-node bifurcation and the dynamics of this ODE is a permanent oscillation. We therefore fix an arbitrary $\muf\in(0,1)$ throughout.

Travelling waves of \eqref{theta-eq} with velocity $c$ solve the ODE, considered in $\mathbb{R}^2$, given by
\begin{equation}
    u_\xi=\phi,\quad \phi_\xi=-c\phi- f(u),\label{e:kink}
\end{equation}
and the fundamental \emph{kinks} are fronts, i.e., heteroclinic connections between the stable states $2\pi (k+1)$, $2\pi k$, $k\in\N$. They are strictly monotone decreasing in $u$ and their unique existence (with $\muf\in(0,1)$) follows, e.g., by phase plane analysis, cf.\ \cite{ErmentroutRinzel}. For $k=0$, we denote the unique translate such that $u(0)=\pi$ by $\varphi$ and note that $c>0$; the basic \emph{antikinks} are translates of the spatial reflection $\varphi(-\xi)$. As a scalar reaction-diffusion equation, fronts in \eqref{theta-eq}, and therefore $\varphi$ as well as $\varphi(-\cdot)$, are orbitally stable \cite{SATTINGER1976312, Fife1977}. We note that due to the periodicity in $u$, there do not exist heteroclinic orbits connecting $2\pi k$ and $2\pi k'$ for $k-k'\neq 1$. 

We use the term kink and antikink more losely for monotone pieces of $u(\cdot,t)$ that connect even multiples of $\pi$. Concerning their positions, we refer to elements in the set $P$ from \eqref{e:pset} as \emph{geometric} positions and will introduce a notion of \emph{analytic} positions in \S\ref{s:positionsdetail}. These types of positions generally differ, but -- as will be shown -- the long term dynamics leads to large distances for which the geometric and analytic positions will be exponentially close to each other. In particular, $P(\varphi(t)) = \{ct\}$ is trivially a single point for all time moving with speed $c$. 

As outlined before, we are interested in the relative motion of sequences kinks and antikinks that resemble superpositions of shifted $\varphi$ and $\varphi(-\cdot)$. We start by considering discontinuous initial data built from kink or antikinks steps,
\[
H^\pm(x) :=  \begin{cases}
2\pi \,, & \pm x < 0\\
\pi\,, & x=0\\
0\, & x>0.
\end{cases}
\]
Notably, the solution with initial data $H^+$ will converge to the above kink $\varphi$, while that with initial $H^-$ converges to the antikink. 
Next, we consider initial positions
\begin{equation}\label{pos_kinkantikink}
    \xi_m^-<\xi_{m-1}^-<\ldots<\xi_1^-<\xi_1^+<\xi_2^+<\ldots<\xi_n^+
\end{equation}
for $m$ kink and $n$ antikink steps that are shifted to these positions via 
\begin{equation*}
H_j^\pm(x):=H^\pm(x-\xi_{j+1}^\pm).
\end{equation*}

For convenience, we smoothen $H^\pm_j$ in such a way that the geometric positions $P(u_0)$ of the resulting $u_0$ coincide with those in \eqref{pos_kinkantikink}. This can be realised by replacing $H_j^\pm$ with a convolution $H_{j,\varepsilon_j^\pm}^\pm:= \rho_{\varepsilon_j^\pm}\ast H_j^\pm$ for a positive, symmetric and smooth mollifier $\rho_\varepsilon$ supported on $[-\varepsilon,\varepsilon]$ and sufficiently small $\varepsilon_j^\pm$ depending on the neighboring initial positions, e.g., $\varepsilon_j^\pm<\frac 1 2 \min\{ | \xi_j^\pm - \xi_{j-1}^\pm|, | \xi_j^\pm - \xi_{j+1}^\pm|\}$ for $2\leq j\leq n-1$ and correspondingly for $j=1,n$. Hence, we set 
\begin{equation}\label{data}
    u_0(x;m,n):=\sum_{j=0}^{m-1}H_{j,\varepsilon_j^-}^-(x)+\sum_{j=0}^{n-1}H_{j,\varepsilon_j^+}^+(x).
\end{equation}
with $n,m\in \overline{\N}_0:=\N\cup\{0,\infty\}$, i.e., possibly infinitely many kink or antikink steps.  Due to the separated smoothened intervals the geometric positions of $u_0$ coincide with \eqref{pos_kinkantikink}, i.e., $P(u_0(x;m,n))=\{\xi_1^-,\ldots,\xi_m^-,\xi_1^+,\ldots,\xi_n^+\}$. 

In the following we omit the dependence on $\varepsilon_j^\pm$ since these do not influence the results. 

\subsection{Well-posedness}

We consider possibly unbounded initial data, which in particular covers the case of infinitely many kinks or antikinks. In order to ensure well-posedness, one option is to consider a weight function $\omega\colon\mathbb{R}\to\mathbb{R}, \omega(x):=C^{-1}e^{-C\lvert x\rvert}$ for some $C>0$, and  
the Banach space
\begin{equation*}
    X_\omega:=\left\{v(\cdot)\in C(\mathbb{R}): \omega v\in L^\infty(\mathbb{R})\right\},\quad\lVert v\rVert_{X_\omega}:=\lVert \omega v\rVert_{\infty}:=\sup_{x\in\mathbb{R}}\lvert v(x)\omega(x)\lvert.
\end{equation*}
\begin{theorem}\label{t:global}
For $f$ from \eqref{theta-eq} the initial value problem
\begin{equation*}
    \begin{cases}u_t=u_{xx}+f(u), & u(x,t)\in\mathbb{R}, (x,t)\in \mathbb{R}\times\mathbb{R}_{>0}\\
    u(x,0)=u_0(x)\in X_\omega, & x\in\mathbb{R}\end{cases}
\end{equation*}
has a unique solution $u\in C^\infty(\mathbb{R}\times\mathbb{R}_{> 0},X_\omega)$. 
\end{theorem}

\begin{proof}
	The proof directly follows nowadays classical techniques, 	and we refer to the milestone monograph \cite{Lady} for further 
details. Let us just briefly recall the main steps.
Rephrasing the initial value problem as an integral equation, the \emph{local} existence of a unique solution can be deduced from a standard fixed point argument. To this end, we consider the operator $\Phi\colon C([0,T],X_\omega)\to C([0,T],X_\omega)$ defined by
\begin{equation}\label{e:varoconst}
    \Phi[u]:=G_t*u_0(x)+J(x,t),\quad J(x,t):=\int_0^t\int_\mathbb{R}G_{t-s}(x-y)f(u(s,y))\, dy\, ds
\end{equation}
where $u(x,\cdot)\in C([0,T],X_\omega)$ for $x\in\mathbb{R}$ and $G(x,t):=G_t(x):=\frac{1}{\sqrt{4\pi t}}e^{-\frac{\lvert x\rvert^2}{4t}},t>0$, is the heat kernel. Moreoever, by proving the H\"older continuity of this solution, a boot-strap argument shows smoothness of the solution both in $x$ and $t$.  

As for ODE with bounded vector fields, the local existence of the solution can be extended to \emph{global} existence ($t>0$) since $T$ is independent of the initial condition.  
\end{proof}

\subsection{Geometric and analytic positions}\label{s:positionsdetail}

Having established global existence, we turn to the specific notions of positions. 
First we note that the geometrically defined set of positions $P(u(x,t;m,n))$ gives locally smooth curves of positions as follows. More details for different types of initial data will be given in the subsequent sections.

\begin{proposition}~
	\label{prop_posfunc}  
	For any $m,n\in \overline{\N}_0$, $n+m>0$, consider the global solution $u(x,t)$ from Theorem~\ref{t:global} with an initial datum $u_0(x;m,n)$ as in \eqref{data}. Then the following holds. For any $t\geq 0$ the set $P(u(t))$ is discrete and, if $n+m<\infty$, consists of at most $n+m$ elements. There is $t_1>0$, $t_1=\infty$ if $nm=0$, such that for $0\leq t<t_1$ the set $P(u(t))$ consist of differentiable curves $ \xi_m^-(t)<\xi_{m-1}^-(t)<\ldots<\xi_1^-(t)<\xi_1^+(t)<\xi_2^+(t)<\ldots<\xi_n^+(t)$ that coincide with the initial positions \eqref{pos_kinkantikink} at $t=0$. 
Moreover, as long as any two such positions are defined, the number of elements in $P(u(t))$ between these cannot increase. 
\end{proposition}

This Proposition in particular yields, at least locally in time, well-defined and regularly varying positions $\xi_j^\pm(t)$;  we refer to $\xi_j^-(t)$ as kink positions and $\xi_j^+(t)$ as antikink positions.

\begin{proof}
This proposition is a consequence of the properties of the number of zeros for linear parabolic
equations, see \cite{RedhefferWalter74,angenent1988zero} for a rigorous exposition 
and we refer to \cite[\S 2.3]{Pauthier_2018} 
for an exposition of the results used next. 

Let $v=\dr_xu$ be the spatial derivative of the solution. Given $u$, it solves a linear parabolic equation, and 
at $t=0$ the minimum of $u$ yields a sign change of $v$ at a locally unique zero. The intervals where $u$ is 
constant occur on monotone parts of $u$ and are thus not associated to non-trivial sign changes of $v$. 

It follows from the zero number principle that there exists 
$T\in(0,\infty]$ such that $v$ has a unique simple zero on $(0,T),$ and has constant sign on 
$(T,\infty).$ Moreover, from parabolic regularity and the implicit function theorem, there 
exists a $C^1$ function $t\mapsto\eta(t)$ such that $v(\eta(t),t)=0,$ for all $t\in(0,T).$
This implies that the set $P(u(t))$ is discrete for all $t\geq 0$ and, again from the implicit function theorem,
that each point lies on a differentiable curve.

Let us now turn to the case $n+m<\infty.$ Let $\beta_0^\pm$ be the value of $u_0$ at $x\to\pm\infty$. 
Then (see \cite[Theorem 4.4.2]{Volpert3},  for instance)  
$$
\beta^\pm(t):=\lim_{x\to\pm\infty}u(x,t).
$$ 
exists for all $t\geq0,$ and are solutions of the initial value problem 
$$
\dot{\beta}^\pm=f(\beta^\pm),\qquad \beta^\pm(0)=\beta_0^\pm
$$
in particular, due to our choice of initial condition (\ref{data}), it follows that these limits 
are constant steady states of the above equation.
Combined with the sign properties of $\dr_xu,$ this proves that $P(t)$ consists 
of at most $n+m$ elements.

It remains to prove that the cardinality of $P(t)$ is non-increasing. We claim that if there exists 
$k_0,t_0$ such that $u(\eta(t_0),t_0)=(2k_0-1)\pi,$ then $u(\cdot,t)>(2k_0-1)\pi$ for all $t>t_0.$ 
This is a consequence of the maximum principle: let $\alpha(t)$ be a solution of the 
initial value problem 
$$
\dot{\alpha}=f(\alpha),\qquad \alpha(t_0)=(2k_0-1)\pi.
$$
Then $t\mapsto\alpha(t)$ is increasing and converges to $2k_0\pi$ as $t\to\infty.$ Let 
$w(x,t)=u(x,t)-\alpha(t).$ Then, taking $u$ as given, $w$ solves a linear parabolic equation on $(t_0,\infty),$  
and $w(\cdot,t_0)\geq0.$ It follows that $w(\cdot,t)>0$ for all $t>t_0,$ which concludes the claim
and the proof.
\end{proof}

In particular, for the solution $u_*(x,t)$ with a single smoothened jump kink initial data $u_*(x,0)=H_\varepsilon$, there is a unique globally defined geometric position $P(u(t)) = \{\xi_*(t)\}$. We denote its speed as $c_*:=\frac{\rmd}{\rmd t} \xi_*$. Recall that the single kink speed is denoted by $c$, and the aforementioned convergence result of $u_*$ to the kink $\varphi$ implies $c_*(t) \to c$ as $t\to\infty$. 

\bigskip
In order to study the evolution of distances between kinks or antikinks in more detail, it is convenient to define the following analytic positions, which however requires sufficiently large initial distances $\xi_i-\xi_{i+1}$ ($1\leq j\leq n$) between neighboring kinks. 

The broader task of deriving laws of motion for localized states in terms of ordinary differential equations ("laws of motion'') dates back at least to the studies on metastable fronts in the Allen-Cahn equation by Carr-Pego and Fusco-Hale, who derived ODEs for the analytic positions of fronts, cf.\ \cite{doi:10.1002/cpa.3160420502,Fusco1989,Fife1977}. This has been explored in various directions, notably to infinitely many metastable pulses in arbitrary dimension \cite{ZelikMielke}; we mainly follow \cite{Ei2002} and \cite{Rossides2015}. 

This allows to derive such ODE rigorously only for sequences of either kinks or antikinks, i.e., monotone initial data, and we therefore restrict attention to \eqref{data} with $m=0$. Since in this case the overall motion of kink initial data is dominated by the drift with velocity $c$, we consider the deviation from this speed by introducing the comoving frame $z=x-c t$, which introduces the term $c\partial_z u$ on the right hand side of \eqref{theta-eq} and yields, in the covering space $\R$, 
\begin{equation}\label{theta-eqz}
    u_t= u_{zz}+c u_z + f(u).
\end{equation}
The corresponding solution $u(z,t)$ is defined globally in $t\geq 0$ by Theorem~\ref{t:global}. We will define analytic positions $\eta_i(t)$ that relate to the geometric positions $\xi_i$ by $ct + \eta_i(t)\approx\xi_i(t)$. 

For sufficiently large initial distances, the analytic positions are defined by writing $u(z,t)$ as 
\begin{equation}\label{ansatz}
u(z,t)=\sum_{i=1}^n\varphi_i(z,t)+w(z,t),\quad \varphi_i(z,t):=\varphi(z-\eta_i(t)),
\end{equation}
with $\eta_i~(1\leqslant j\leqslant n)$ uniquely defined through the following orthogonality condition on $w$ as detailed in Appendix~\ref{derivationODE}.
Let $L_i:=\partial_z ^2+c\partial_z+f'(\varphi_i)$ denote the linearized operator of the right hand side of \eqref{theta-eqz} in $\varphi_i$, and $L_i^*:=\partial_z^2-c\partial_z+f'(\varphi_i)$ its adjoint. 
The remainder term $w$ is now supposed to be orthogonal to the adjoint eigenfunctions $e_i^*:=e^{c(z-\eta_i)}\varphi_i'$, i.e.
\begin{equation}\label{orthcond}
\langle w,e_i^*\rangle=\int_\mathbb{R} e^{c(z-\eta_i)}w(z,t)\varphi_i'(z,t)\, dz=0,\quad 1\leqslant i\leqslant n.
\end{equation}
For initial data as in \eqref{data}, $w$ is initially nonzero, so that for the study of analytic positions it is natural -- though not necessary -- to replace $H^\pm_{k,\varepsilon}$ by the fronts, i.e., $\varphi_i$ as in \eqref{ansatz}. Then $w(z,0)\equiv0$ and thus remains small at least for short time. In order to control $w$ also for larger times, we assume that the kinks are initially well-separated, i.e., $\lvert \eta_i(0)-\eta_{i-1}(0)\rvert\gg 1~(2\leqslant i\leqslant n)$.

\begin{remark}\label{r:kinkpos}
For initial data $v_0$ as in \eqref{ansatz}, the set $P(v_0)$ does not coincide with the initial positions \eqref{pos_kinkantikink}, though the error is exponentially small in the minimal distance. Moreover, for $nm>0$, the set $P(v_0(\cdot;m,n))$ may contain spurious points such that one generally needs to assume a minimal distance between $\xi_1^-$ and $\xi_1^+$.
\end{remark}

\begin{remark}
For any fixed $t$ and $i$, the analytic kink position $\eta_i(t)$ equals the (shifted) geometric position $\xi_i(t)-c t$ if and only if $w(\xi_i(t),t)=(n-1)\pi - \sum_{j=1,j\neq i}^n\varphi(\xi_j(t)-ct-\eta_j(t))$, as can be seen from \eqref{ansatz}. In particular, all geometric and analytic positions coincide if and only if $w$ vanishes simultaneously at all geometric positions. 
However, kinks interact eventually repulsively, i.e., their distances eventually increase (see \S\ref{boundedmonotone} below for details) and this implies $\lVert w(\cdot,t)\rVert\to 0$ as $t\to\infty$ (in $L^2$ or $L^\infty)$ so that both definitions asymptotically coincide. 
\end{remark}

\section{Bounded monotone initial kink data}\label{boundedmonotone}
In this section we consider bounded monotone data which are composed of kinks and analyse the distances between neighboring kinks in terms of the geometric as well as analytic positions. 
We note that by spatial reflection the discussion equally applies to bounded sequences of antikinks. First we track geometric positions via the comparison principle, which applies for any initial distance, but only constrains the positions to lie within certain intervals, referred to as \textit{gaps}, that also depend on the initial data. The analytic positions provide more specific laws of motion that apply immediately for initial data with sufficiently distant positions, or -- more abstractly -- from some point in time onward with a distribution of positions for which we just know the positions up to the gaps. 

The latter relies on the results by Pol\'{a}\v{c}ik and Risler, which state that front-like initial data converge to a terrace whose speeds converge and whose distances \emph{eventually} diverge, albeit without a quantitative estimate. For our purposes, this can be summarized as follows.

\begin{theorem}\label{t:pol}
\label{t:polacikrisler} cf.~\cite{polacik2016propagating,Risler}\\
Let 
$-\pi< u_0 < 2\pi (k+1)~(k\in\mathbb{N})$ be an initial datum with $\lim_{x\to-\infty}u_0(x)=2\pi k$, $\lim_{x\to\infty}u_0(x)=0$ and corresponding solution $u$. Then there exist $C^1$ functions $\xi_1,\ldots,\xi_k$ on $\mathbb{R}$ satisfying 
\begin{equation*}
 \lim_{t\to\infty}\xi_j'(t)=0~(j=1,2,\ldots,k), \qquad  \lvert\xi_j(t)-\xi_{j+1}(t)\rvert\to\infty~(j=1,2,\ldots,k-1) 
\end{equation*}
such that the solution $u(\cdot,t)$ converges to the corresponding terrace:
\begin{equation*}
\lV u(\cdot,t)-\sum_{j=1}^k\vp(\cdot-ct-\xi_j(t))\rV_\infty\underset{t\to\infty}{\longrightarrow}0
\end{equation*}
\end{theorem}

\begin{figure}[H]
    \centering
    \begin{tabular}{cc}
   \hspace{-1.3cm} \includegraphics[scale=.7]{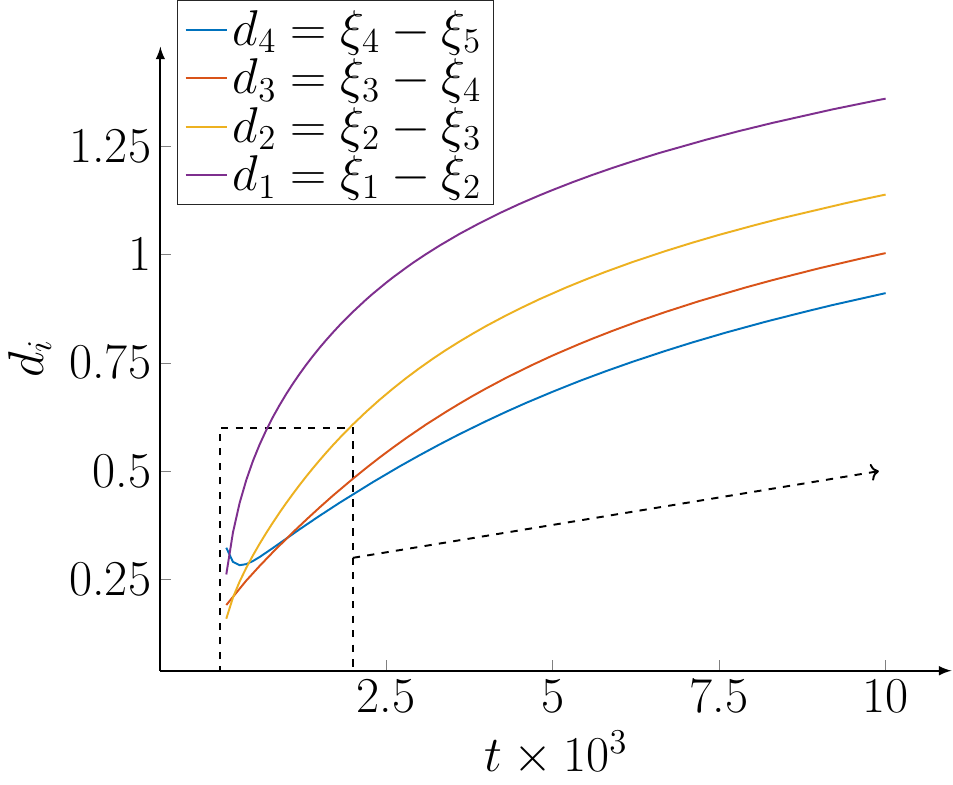} & \hspace{-0.3cm}\includegraphics[scale=.7]{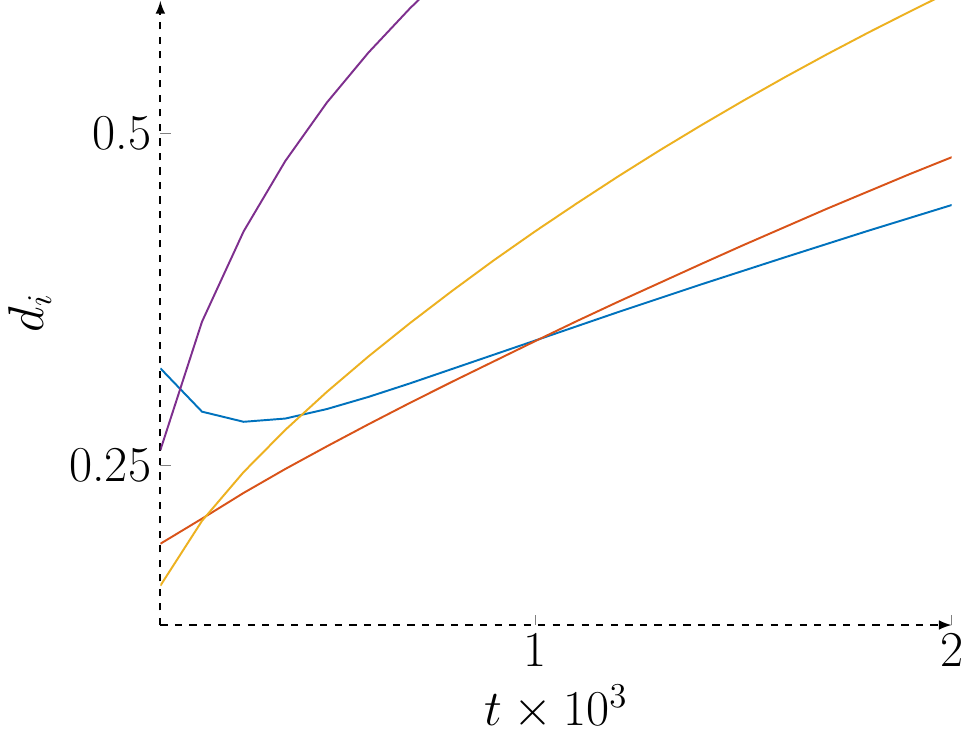}\\
   (a) & (b)\\
   \hspace{-1.3cm} \includegraphics[scale=.7]{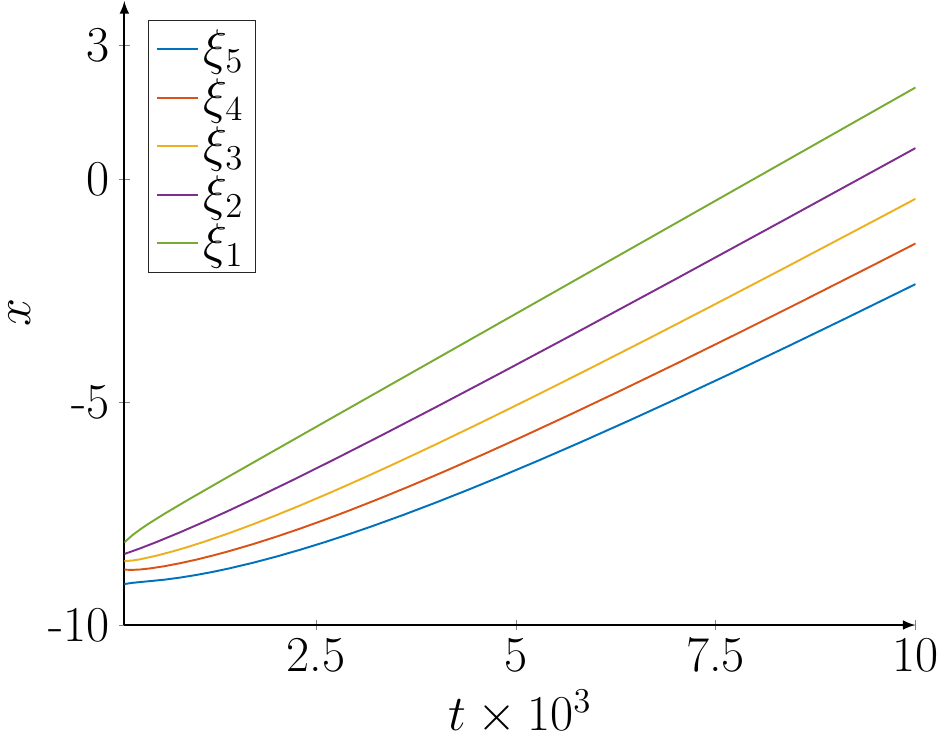} & \hspace{-0.3cm}\includegraphics[scale=.7]{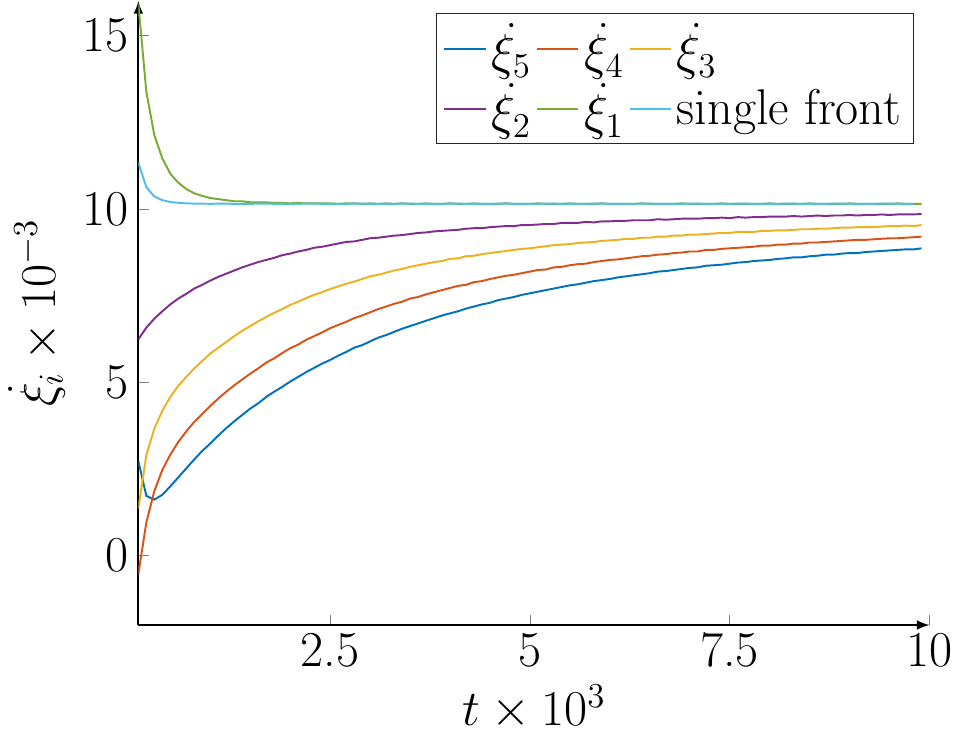}\\
   (c) & (d)
    \end{tabular}
   \caption{Simulation of kink distances (a) (with zoom (b)), positions (c) and speeds (d), starting from an initial condition with five kinks and $\muf=0.2$. As can be seen from (a) and (b), the distances need not be monotone functions, but are eventually ordered, cf. \S\ref{s:monotone_ODE}. Moreover, the speeds of the kinks converge to the single front speed (d).}
   \label{fig:ordering}
\end{figure}

\begin{remark}
In Prop.~\ref{prop_non-equi} below we give a lower bound on the distances which is of course far from optimal in the asymptotics $t\to\infty$. Naively, one might suppose that all distances monotonically increase; however, this is not true in general, as the results in \S\ref{s:monotone_ODE} show, and numerical simulations illustrate, cf. Fig.~\ref{fig:ordering}~(a) and (b). 
\end{remark}

In order to describe the long-time behaviour of bounded solutions under kink-antikink annihilations, we consider limit sets. For bounded monotone solutions, this has already been done -- for a much broader setup -- in \cite{polacik2016propagating} from which we take the following definition of the limit set
\begin{equation*}
    \Omega(u):=\{v: u(\cdot+x_n,t_n)\to v~\textrm{ for some sequences }t_n\to\infty\textrm{ and }x_n\in\R\},
\end{equation*}
where $u\in L^\infty(\mathbb{R}\times\mathbb{R}^+)$ and the convergence is in $L^\infty_{loc}(\mathbb{R})$ (locally uniform convergence). Compared to the standard definition of $\omega$-limit set, this allows to observe any finite piece of the graph of $u(\cdot,t_n)$. For the special situation of Theorem \ref{t:polacikrisler}, the set is given by
\begin{equation}\label{Omegaset}
    \Omega(u)=\{\varphi(\cdot-\tau): \tau\in\R\}\cup\{2\pi j: j=0,1,\ldots,k\},
\end{equation}
cf. \cite{polacik2016propagating}. Our results add the description of $\omega$ and $\Omega$ in case $m+n<\infty$ and $mn\neq 0,$ even though they are consequences of \eqref{Omegaset} after pairwise annihilations of kinks and antikinks.

\subsection{Qualitative aspects: comparison principle}\label{subsection2.1}
Let $u_0(x;n):=u_0(x;n,0)$ be a monotone initial datum \eqref{data} without antikinks and
associated global solution $u(x,t;n)$; recall the positions are globally defined and differentiable according to Proposition~\ref{prop_posfunc}. Since the data has kinks only we omit the superindex $+$ and also introduce the following notation for the speeds of positions, the nearest distances and the minimal distance (starting from a given position):
\begin{align}
   c_i(t)&:=\frac{d}{dt}\xi_{i}(t),\quad 1\leqslant i\leqslant n,\label{speed}
    \tag{\theequation$$}
    \stepcounter{equation}\\
    d_j(t)&:=\xi_j(t)-\xi_{j+1}(t),\quad 1\leqslant j\leqslant n-1\label{distance}
    \tag{\theequation$$}
    \stepcounter{equation}\\
    \underline{d}_i(t)&:=\min_{i\leqslant j\leqslant n-1}d_j(t),
\end{align}
as well as $d_{\min}(t):=\underline{d}_1(t)$.

\begin{figure}
    \centering
    \includegraphics[scale=0.9]{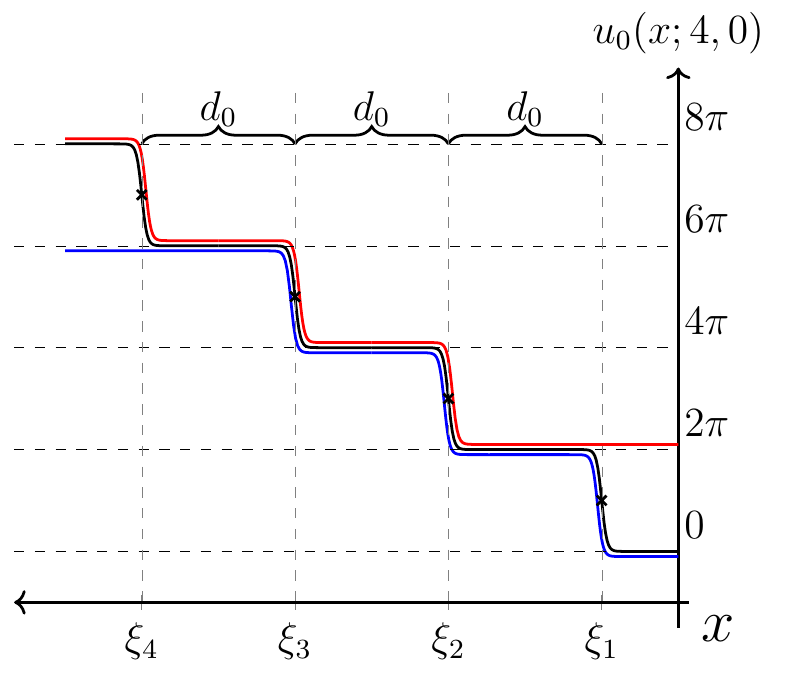}
    \caption{Sketch of initial datum $u_0(x;4)$ (black) with four antikinks at equidistant positions $\xi_{i}(0), i=1,2,3,4$, such that $u(\xi_{i}(0),0)= (2 j + 1)\pi $. 
    The blue initial datum is given by $\underline{u}_0$ with position $\tau_{j}(0)=\xi_{j}(0)$ for $j=1,2,3$ while the red one is given by $\bar{u}_0$ with positions $\psi_{j}(0)=\xi_{j+1}(0)$ for $j=1,2,3$. For the purpose of illustration, the curves are plotted slightly below and above $u_0(x;4)$, respectively. 
} 
\label{fig:decreasing-staircase}
\end{figure}

\textbf{Equidistant kinks.}
Let us first consider initial data composed of equidistant antikinks, cf.\ Figure~\ref{fig:decreasing-staircase}, as it is straightforward to construct sub- and supersolutions in order to compare the speeds {\let\pm=-\eqref{speed}} and, consequently, to find a uniform lower bound on the distances {\let\pm=-\eqref{distance}}.
\begin{proposition}\label{prop_equi}
Let $u(x,t;n)~(n\geqslant 2)$ be the solution with an initial datum $u_0(x;n)$, where $d_i(0)=d_0$ for all $1\leqslant i\leqslant n-1$ for some constant $d_0>0$. Then $c_1(t)\geqslant c_*(t)\geqslant c_n(t)$ for all $t\geqslant 0$, and $c_1(0)\geqslant c_2(0)\geqslant\ldots\geqslant c_n(0)$ . In particular, initially all distances $d_i(0)$ are non-decreasing in $t$ and $d_i(t)\geqslant d_0$ for all $t>0$.
\end{proposition}

\begin{proof}
We will consider solutions $u(x,t;k)$ with initial $u_0(x;k)$ for the same positions $\xi_j$, $j=1,\ldots,k$, but $k\leq n$, and compare the speeds of kinks which we therefore denote as $c_j^k$, where we omit the argument $t$ for readability; hence, $c_* = c_1^1$. 

We first prove the statement for $n=2$. In this case, the initial datum $u_0(x;n)$ is sandwiched by $H(x-\xi_1)\leqslant u_0(x;n)\leqslant H(x-\xi_2)+2\pi$ for all $x\in\mathbb{R}$. By the comparison principle, for the associated solution $u_*$ we have $u_*(x-\xi+1,t)\leqslant u(x,t;n)\leqslant u_*(x-\xi_2,t)+2\pi$ for all $(x,t)\in\mathbb{R}\times\mathbb{R}_+$. Since both $ H_0$ and $ H_1$ move with the same speed $c_*$, the speeds of the kinks satisfy $c_1^{2}\geqslant c_*\geqslant c_2^{2}$ for all $t\geqslant 0$. 

For $n\geqslant 3$, we choose sub- and supersolutions composed of $n-1$ antikinks (cf. Figure~\ref{fig:decreasing-staircase}). More precisely, let $u_0(x;n-1)$ be the initial condition composed of $n-1$ kinks with the same $n-1$ positions, $\xi_{i}(0)$ for $2\leqslant i\leqslant n$, as $u_0(x;n)$. The comparison solutions arise from the initial data
\begin{align*}
&\underline{u}_0(x):=u_0(x-d_0;n-1), &
&\bar{u}_0(x):=u_0(x;n-1)+2\pi.
\end{align*}
In particular, the kink positions of $\underline{u}_0$ are $\xi_{i}(0)$ for $1\leqslant i\leqslant n-1$ and those of $\bar{u}_0$ are  $\xi_{i}(0)$ for $2\leqslant i\leqslant n$. 

Since $\underline{u}$ and $\bar{u}$ are translates of each other, their positions $\underline{\xi}_i$, $\bar{\xi}_i$, $1\leqslant i\leqslant n-1$ are equal up to shift by $d_0$ and have the same speeds $c_i^{n-1}$, for all $t\geqslant 0$. By construction, $\underline{u}_0\leqslant u_0\leqslant \bar{u}_0$ and thus the corresponding solutions satisfy $\underline{u}\leqslant u(\cdot;n)\leqslant \bar{u}$ for all $(x,t)\in\mathbb{R}\times\mathbb{R}_+$ by the comparison principle. We can therefore compare the speeds $c_j^{n}$  with the speeds $c_i^{n-1}$ of positions in $\underline{u}$ and $\bar{u}$, respectively, which leads to the following relations for all $t\geqslant 0$:
\begin{align*}
    &c_1^{n}\geqslant c_1^{n-1}, 
    &&c_j^{n-1}\leqslant c_j^{n}\leqslant c_{j-1}^{n-1}\textrm{ for }1<j<n, 
    &&c_n^{n}\leqslant c_{n-1}^{n-1}.
\end{align*}
Hence, $c_i^{n}\geqslant c_i^{n-1}\geqslant c_{i+1}^{n}$ for $i=1,\ldots,n-1$. By iterating this construction of sub- and supersolutions, $c_1^{n}\geqslant c_1^{n-1}\geqslant\ldots\geqslant c_1^{2}\geqslant c$ and, analogously, $c_n^{n}\leqslant c_{n-1}^{n-1}\leqslant\ldots\leqslant c_2^{2}\leqslant c$ for all $t\geqslant 0$.
\end{proof}

We note that the proof of Prop.~\ref{prop_equi} shows that there is a hierarchy of speeds when removing kinks on the left or right. Up to equality in the bounds, kinks for equidistant data spread out, and the more  kinks there are, the faster the spreading can be.

\begin{figure}[H]
    \centering
    \includegraphics[scale=0.9]{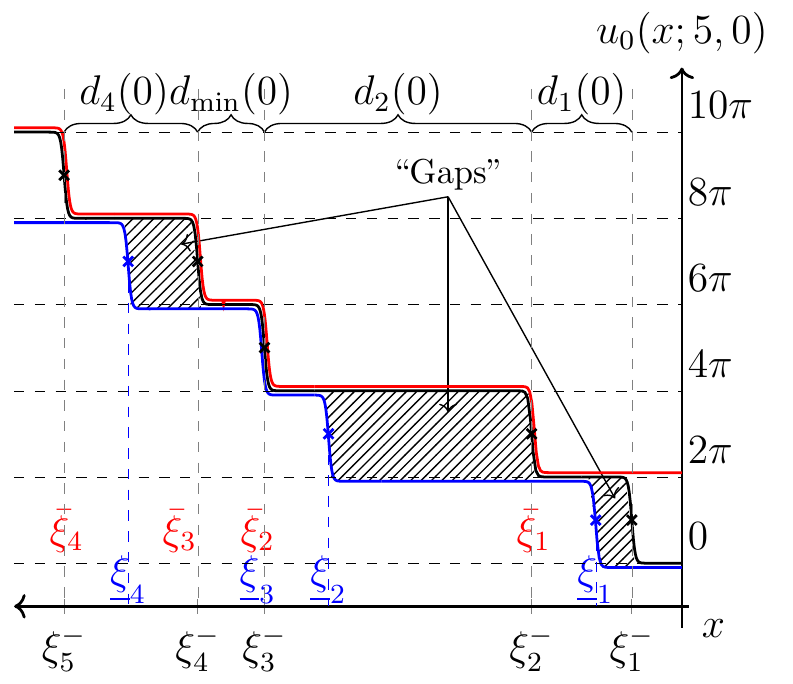}
    \caption{Sketch of gaps (hatched): $\underline{\xi}_i(0)<\xi_i(0)$ for $i=1,2,4$. The associated speeds need not be ordered. In particular, the distances may be increasing or decreasing. }\label{fig:decreasing-non-equid}
\end{figure}

\textbf{Non-equidistant kinks.} 
For more general initial data, only the overall distance $\xi_n(t)-\xi_1(t)$ and the smallest distance $d_{\min}(t)$ can be controlled by our approach of sub- and supersolutions: lower bounds can be inferred for the distances only up to ``gaps'' rooted in the initial data. In these gaps the ordering of the associated positions and speeds cannot be further constrained by this method, cf. Fig.~\ref{fig:decreasing-non-equid}. 
In view of the upcoming analysis based on analytic positions, this is not surprising since distances may behave non-monotonically. 

 \begin{proposition}\label{prop_non-equi}
The solution to any initial $u_0(x;n)$, $n\geqslant 2$, satisfies the following for all $t\geqslant 0$:
\begin{align}
    &\frac{\rmd}{\rmd t}(\xi_n(t)-\xi_1(t))\geqslant 0\,,\quad c_1(t) \geqslant c_*(t) \geqslant c_n(t),\nonumber\\
    & d_i(t)\geqslant d_{\mathrm{min}}(0)\quad\forall t>0, \;1\leqslant i\leqslant n-1.\label{lowerbound}
\end{align}
Moreover, if $d_{\min}(0) = d_j(0)$ for some $1\leqslant j\leqslant n-1$ and $t\geqslant 0$ then $d_j(t)$ is non-decreasing for all $t\geqslant 0$. 
In particular, $d_{\min}(t)$ is non-decreasing at least as long as the minimal distance is realised at the same index as initially. 
\end{proposition}
\begin{proof}
The proof is analogous to that of Prop.~\ref{prop_equi}. Comparing $u$ with single kinks shifted to $\xi_1$ as a subsolution and to $\xi_n$ as a supersolution we immediately infer $c_1(t) \geqslant c_*(t) \geqslant c_n(t)$ so that $\frac{\rmd}{\rmd t}(\xi_n(t)-\xi_1(t))\geqslant 0$. 
Next, we replace $d_0$ in the proof of of Prop.~\ref{prop_equi} by $d_{\min}(0)$ and use again $\bar{u}_0(x):=u_0(x-d_{\min}(0);n-1)$, as well as $\underline{u}_0:=u_0(x;n-1) +2\pi$, cf.\ Fig.~\ref{fig:decreasing-non-equid}. By the comparison principle this implies for all $t\geqslant 0$ the relations $c_i^n\leqslant c_{i-1}^{n-1}$ for $2\leqslant i\leqslant n$ from the supersolution, but the subsolution only yields for $j$ such that $d_{\min}(0)=d_j(0)$ the relation $c_j^{n-1}\leq c_j^n$. Taken together we obtain, for all $t\geqslant 0$,
\[
c_j^n\geqslant c_j^{n-1}\geqslant c_{j+1}^n,
\]
and therefore $\frac{\rmd}{\rmd t}d_j(t)\geq 0$. In particular, $\frac{\rmd}{\rmd t}d_{\min}(t)\geqslant 0$ as long as the minimal distance is between $\xi_j(t)$ and $\xi_{j+1}(t)$, but at least for $t=0$. 
\end{proof}

As this use of sub- and supersolutions is limited by the described gaps, a substantial part of the next section is devoted to a deeper analysis of the behaviour of the distances in terms of the analytic positions.

\subsection{Quantitative aspects: projection scheme}\label{s:monotone_ODE}

The above discussion, and in particular Proposition \ref{prop_non-equi}, gives a first qualitative overview of 
interactions between stacked fronts traveling at the same asymptotic speed $c$.
The purpose of this section is to obtain deeper insight into the relative laws of motion of the front positions,  
which is however valid for sufficiently large distances only.
As mentioned in \S\ref{s:positionsdetail}, we follow an approach that has developed out of the analysis of 
meta-stability of fronts in Allen-Cahn equation in \cite{doi:10.1002/cpa.3160420502,Fusco1989}. 

Recall the notion of analytic positions from \eqref{ansatz}, \eqref{orthcond}, and denote $\eta:=(\eta_1,\ldots, \eta_n)$ as well as
\begin{equation*}
\delta_j:=\eta_{j}-\eta_{j+1}\;,\qquad\delta_{\min}=\min\{\delta_j,1\leq j\leq n-1\}.
\end{equation*}

Initial data \eqref{ansatz} with $w\equiv0$ form the nested  $n$-dimensional manifolds (with $\delta=\delta_{\min}$)
\[
\calK_\delta:= \left \{ \sum_{j=0}^n \vp(\cdot+z_0-\eta_j) : \eta_1>\eta_2>\cdots>\eta_n, \; \eta_i-\eta_{i+1} > \delta, i=1,\ldots, n; \, z_0\in\R \right\}.
\]

Theorem~\ref{t:pol} shows that kinks interact eventually repulsively, which can be understood as asymptotic stability of the boundary   of $\calK_0$. In this section we will show that, for any sufficiently large $\delta$, $\calK_\delta$ parametrises an $n$-dimensional invariant manifold $\calM$, and we study the reduced dynamics on it. Moreover, $\calM$ is globally attracting for initial data as in Theorem~\ref{t:pol}. In particular, any initial datum with initial data from geometric positions \eqref{data} is contained in its basin of attraction, if the initial geometric distances are sufficiently large. 

This analysis relies on the following results, which can be inferred from \cite{Ei2002}. Although the latter is not intended to specifically study stacked fronts as in his paper, but rather pulses, our situation can be dealt with by exactly the same methods. 

We reformulate the results that we use in our notation and towards our purposes, and refer in particular to Theorems 2.1, 2.3 and 4.1 in \cite{Ei2002} for proofs. 
We also provide in Appendix \ref{derivationODE} some more details and a derivation of the reduced ODE system, in order to close the gap to \cite{Ei2002}. 

The basic observation is the following well known decomposition in such neighborhoods, which readily follows from an implicit function theorem, cf.\ Appendix \ref{derivationODE}. Here we consider neighbourhoods of $\calK_\delta$ defined as $\calU_{\eps,\delta}:=\calK_\delta+\left\{v\in H^2(\R) : \|v\|_{H^1}<\eps \right\}$.

\begin{lemma}\label{lem:decomp}
There exist $\delta^*, \eps^*>0$ such that for any $v\in \calU_{\eps^*,\delta^*}$, there exist unique $\eta_1,\cdots \eta_n\in \R^n$ and a unique function $w(z)\in H^2(\R)$ satisfying
\begin{equation}\label{eq:17x}
v(z)=\sum_{j=0}^{n}\vp(z-\eta_j)+w(z)\;,\quad
\int_\mathbb{R}e^{cz}w(z)\vp(z-\eta_j)\rmd z=0,\textrm{ for all }1\leq j\leq n.
\end{equation}
\end{lemma}

We note the following direct consequence of Lemma~\ref{lem:decomp} together with well-posedness.

\begin{proposition}\label{prop:Ei1}
For any solution $u(z,t)$ of (\ref{theta-eqz}) whose initial data lies in $\calU_{\eps^*,\delta^*}$ we have $T_*:=\sup\{ t : u(z,s) \in \calU_{\eps^*,\delta^*}, 0\leq s\leq t\}>0$ and there exist unique functions $\eta_1(t),\cdots \eta_n(t)\in C^1(0,T_*),$ and  
$w(z,t)\in C^1((0,T_*)\times\mathbb{R})$ such that \eqref{ansatz} and \eqref{orthcond} hold for all $t\in(0,T_*)$. 
\end{proposition}

In particular, initial conditions $u_0(x;0,n)$ from (\ref{data}) with $d_{\min}$ sufficiently large lie in $\calU_{\eps^*,\delta^*}$ so that Proposition \ref{prop:Ei1} applies for those as well.

Building on the decomposition \eqref{ansatz}, the following proposition, which is a combination of results from \cite{Ei2002}, cf.\ Appendix \ref{derivationODE}, gives a quantitative description of the relative motions as a reduced ODE for the analytic positions $\eta_j.$ Relevant for the laws of motion are in particular the eigenvalues 
of the linearization of \eqref{e:kink} in the asymptotic states $(\vp,\vp')=(0,0)$ and $(\vp,\vp')=(2\pi,0),$ 
which we readily determine as
\begin{equation}\label{eq:42}
\lambda:=\frac{-c-\sqrt{c^2-4f'(0)}}{2}<0,\ \mu:=\frac{-c+\sqrt{c^2-4f'(0)}}{2}>0;
\end{equation}
recall $f'(2\pi)=f'(0)$ by periodicity of $f,$ and note $\lambda=-\mu-c.$

\begin{proposition}\label{prop:Ei2}
There exist $\delta_*\geq \delta^*, \eps_*\leq \eps^*$ such that $\calU_{\eps_*,\delta_*}$ contains an $n$-dimensional exponentially attractive locally invariant manifold $\calM$, which is a graph over $\calK_{\delta_*}$ in terms of corrections $w=w(\eta)$ satisfying \eqref{eq:17x}. 
The reduced $n$-dimensional dynamics of \eqref{theta-eqz} on $\calM$ is given by the following ordinary differential equation system, which is defined for $\delta_{\min}\geq \delta_*$ and where $\alpha, a_{\textrm{R}},a_{\textrm{L}}>0,$ and $\mathcal{R}_j(\eta)$, $(1\leq j\leq n)$ are $C^1-$functions; furthermore $\mu+\alpha>-\lambda$ for $\muf\approx 0$. 
\begin{equation}\label{eq:20}
\left\{
\begin{array}{llr}
\eta_1'(t)= & a_{\textrm{L}}e^{\lambda(\eta_1-\eta_2)}+\mathcal{R}_1(\eta) & \\
\eta_j'(t) = & a_{\textrm{L}}e^{\lambda(\eta_j-\eta_{j+1})}-a_{\textrm{R}}e^{-\mu(\eta_{j-1}-\eta_j)}+\mathcal{R}_j(\eta), & 2\leq j\leq n-1, \\
\eta_n'(t) = & -a_{\textrm{R}}e^{-\mu(\eta_{n-1}-\eta_n)}+\mathcal{R}_n(\eta). & 
\end{array}\right.
\end{equation}
The remainder terms $\mathcal{R}_j(\eta)$ and the correction term $w(\eta)$ satisfy, for some constant $C>0$,
\begin{equation}\label{eq:21}
|\mathcal{R}_j|,\lV w(\eta)\rV_{H^2(\R)}\leq C e^{-(\mu+\alpha)\delta_{\min}}.
\end{equation}
\end{proposition}
\begin{remark}\label{rem:orders}
Remark that we do not claim that $\mu+\alpha\geq\lambda$ for all $\muf\in(0,1)$, which means that the remainder terms $\mathcal{R}_j$ and $w$ are not necessarily higher order compared with the $a_{\textrm{L}}$-terms whose exponential rate is $\lambda$. 
 \end{remark}

In the remainder of this section we analyse \eqref{eq:20} in more detail and in particular infer that, if the initial distance are large enough, then the validity constraint $\delta_{\min}\geq \delta_*$ is satisfied for all $t\geq 0$.  This is of course consistent with Theorem~\ref{t:pol}, which moreover implies that all solutions with certain initial data eventually satisfy \eqref{ansatz}, \eqref{orthcond} and obey \eqref{eq:20}.

\begin{theorem}\label{t:valid}
There is $\delta_1\geq \delta_*$ such that $\calU_{\eps_*,\delta_1}$ forward invariant for \eqref{theta-eqz}. In particular, the decomposition \eqref{ansatz}, \eqref{orthcond} and the ODE \eqref{eq:20} are valid for all $t\geq 0$.
Moreover, for any solution $u(z,t)$ with initial datum as in Theorem~\ref{t:pol}, there is $t_1\geq 0$ such that $u(\cdot, t_1)\in \calU_{\eps_*,\delta_1}$.  
\end{theorem}
\begin{proof}
The proof is given in \S\ref{s:perturbed}.
\end{proof}

The analysis of \eqref{eq:20} is facilitated by normalizing $\td_i:=\mu\delta_i(t/\mu)$ and setting $\td_{\min}:=\min\{\td_1,\ldots,\td_n\}$, $\varepsilon:=\min\{\alpha,c\}/\mu$, so that \eqref{eq:20} takes the form
\begin{equation}\label{perturbed}
\begin{aligned}
    \td_1'&=e^{-\td_1}+g_{1}e^{-(1+\varepsilon)\td_{\min}},\\
    \td_j'&=e^{-\td_j}-e^{-\td_{j-1}}+g_{j}e^{-(1+\varepsilon)\td_{\min}},
\end{aligned}
\end{equation}
where $g_j$, $1\leq j\leq n$, are bounded. Note that while $\td_{\min}$ is in general only Lipschitz continuous due to the minimum function, the products $g_{j}e^{-(1+\varepsilon)\td_{\min}}$, $j=1,2$ are smooth. For convenience, we extend their domain beyond the definition limit $\td_{\min} \geq \delta_*$ to smooth functions on $\R^n_+$ with $g_j$ bounded.

Our aim is to analyze the temporal ordering and asymptotics of the distances $\td_j, j=1,2\ldots,n$. To this end, we consider system \eqref{perturbed} as a perturbation of the system with $g_j\equiv 0$ for all $1\leq j\leq n$, which we refer to as the \emph{unperturbed} system. 

\subsubsection{The unperturbed distance system} 
For the unperturbed system, i.e. the system, 
\begin{equation}\label{e:unperturbed}
\begin{aligned}
    \td_1'&=e^{-\td_1},\\
    \td_j'&=e^{-\td_j}-e^{-\td_{j-1}},\quad 2\leqslant j\leqslant n,
    \end{aligned}
\end{equation}
we directly get $\td_1'>0$, that is the rightmost distances is increasing. Moreoever, the distance $\td_j~(2\leqslant j\leqslant n)$ increases if (and only if) $\td_j<\td_{j-1}$, i.e., if the distance between the $j$-th kink and its left neighbor kink is smaller than the distance to its right neighbor. In particular, the smallest distance at time $t\geqslant 0$ increases. This motivates the question whether the distances are eventually ordered and which asymptotics are implied by this ordering.
\begin{theorem}\label{thm:orderingunperturbed}
There exists $T>0$ such that $\td_n(t)<\td_{n-1}(t)<\ldots <\td_1(t)$ for all $t\geqslant T$.
Moreover, $\lim_{t\to\infty}\td_k(t)=\infty$ for all $1\leqslant k\leqslant n$. 
\end{theorem}
\begin{proof}
Let us first prove the divergence of the distances by induction over $n$. The statement is clear for $k=1$ since $\td_1(t)=\log(t+C), C>0$. Suppose the statement holds for $\delta_k$ with $2\leqslant k\leqslant n$ and, by contradiction, $\td_{n+1}\leqslant M$ for some $M>0$. Then $\td_{n+1}'(t)\geqslant e^{-M}-e^{-\td_n(t)}>0$ for $t$ large enough. Consequently, $\td_{n+1}$ converges, hence $\lim_{t\to\infty}\td_{n+1}'(t)=0$ and thus $\td_n(t)\to M$ as $t\to\infty$, contradicting the induction hypothesis.

As to the eventual ordering, we first show the existence of some $T_1>0$ such that $\td_2(t)<\td_1(t)$ for all $t>T_1$. To this end, assume by contradiction that for all $t>0$ there exists some $t_1>t$ with $\td_2(t_1)\geqslant \td_1(t_1)$, i.e. $\td_2'(t_1)\leqslant 0$. By the divergence of the distances, there exists $t_2>t_1$ with $\td_2'(t_2)>0$. Since $\td_1$ increases monotonically, this implies that $\td_2$ intersects the $\td_2$-nullcline $\{(\td_1,\td_2): \td_1=\td_2\}$ infinitely often. However, this is impossible since the set $\{(\td_2,\td_1): \td_2<\td_1\}=\{(\td_2,\td_1): \td_2'>0\}$ 
is forward invariant. In particular, if $\td_2(t)=\td_1(t)$ for some $t\geqslant 0$, then $\td_2(t')<\td_1(t')$ for all $t'>t$. Analogously, there exists $T_2>T_1$ with $\td_3(t)<\td_2(t)$ for all $t>T_2$. Likewise, for all remaining pairs of distances, there exists some suitable time $T_i>T_{i-1}$. Setting $T:=T_1$ concludes the proof.    
\end{proof}

In order to get a result similar to Theorem~\ref{thm:orderingunperturbed} for the perturbed system, we first reformulate the problem by applying a polar blow-up transformation; this enables us to apply a perturbation argument.

\subsubsection{Polar blow-up }
Substituting $z_j=e^{-\td_j}$ into \eqref{perturbed}, and setting $z_{\max}:=e^{-\td_{\min}}$ we obtain the system 
\begin{align*}
    z_1'&=-z_1^2-g_{1}z_1 z_{\max}^{1+\varepsilon},\\
    z_j'&=-z_j^2+z_{j-1}z_j-g_{j}z_j z_{\max}^{1+\varepsilon},\quad 2\leqslant j\leqslant n
\end{align*}
which posssesses the non-hyperbolic equilibrium $\vec{0}=(0,0,\ldots,0)^\intercal\in\mathbb{R}^n$. The polar blow-up transformation $\vec{z}=(z_1,z_2,\ldots,z_n)^T=r(t)\vec{\Psi}(t)$ with $r\in\mathbb{R}$ and  $\vec{\Psi}=(\Psi_1,\Psi_2,\ldots,\Psi_n)^T\in S^{n-1}$ yields
\begin{equation}\label{e:zeq}
\begin{aligned}
    z_1'&=-r^2\Psi_1^2-r^{2+\varepsilon}g_{1}\Psi_1\Psi_{\max}^{1+\varepsilon},\\
    z_j'&=-r^2\Psi_j(\Psi_j-\Psi_{j-1})-r^{2+\varepsilon}g_{j}\Psi_j\Psi_{\max}^{1+\varepsilon},\quad 2\leqslant j\leqslant n,
\end{aligned}
\end{equation}
where $z_{\max} = r \Psi_{\max}$, i.e., for each $t$ there is $j$ such that  $\Psi_{\max}(t)=\Psi_j(t)$.
Note that the equilibrium $\vec{0}$ corresponds to $r=0$. The inner product with $\vec{\Psi}$ gives
\begin{equation}\label{radialperturbed}
    \langle{\vec{z}\,'},\vec{\Psi}\rangle=r'(t)=\sum_{k=1}^nz_k'\Psi_k=-r^2\Sigma_n-\calO(r^{2+\varepsilon})
\end{equation}
where $\Sigma_n:=\Psi_1^3-\sum_{k=2}^n\Psi_k^2(\Psi_{k-1}-\Psi_k)$. 
Using $z_k'=r'\Psi_k+r\Psi_k'$ together with \eqref{radialperturbed} in \eqref{e:zeq} we may divide left and right hand sides by $r$ to obtain
\begin{align*}
  \Psi_1'&=r \Psi_1(\Sigma_n-\Psi_1) +\calO(r^{1+\varepsilon}),\\
\Psi_j'&=r \Psi_j(\Psi_{j-1}-\Psi_j+\Sigma_n) +\calO(r^{1+\varepsilon}).
\end{align*}
Dividing the right hand side by $r$ results in the desingularised radial and angular equations
\begin{align}
    r'&=-\Sigma_n r-\calO(r^{1+\varepsilon}),\label{radialperturbeddyn}\\
    \Psi_1'&=\Psi_1(\Sigma_n-\Psi_1)+\calO(r^\varepsilon),\label{angularperturbeddyn1}\\
    \Psi_j'&= \Psi_j(\Psi_{j-1}-\Psi_j+\Sigma_n) +\calO(r^\varepsilon)\label{angularperturbeddyn2},
\end{align}
which, for $r>0$, has the same trajectories as \eqref{radialperturbed} and solutions are related by time rescaling.

Analogous to \eqref{e:unperturbed}, we consider system \eqref{radialperturbeddyn}-\eqref{angularperturbeddyn2} as a perturbation of the \emph{unperturbed polar system} with $g_j$ set to zero ($1\leq j\leq n$), i.e., zero error terms in \eqref{radialperturbeddyn}-\eqref{angularperturbeddyn2}, which means
\begin{align}
  r'&=-\Sigma_n r,  \label{rdynamics}\\
    \Psi_j'&=\Psi_j(\Sigma_n-\Psi_j+\Psi_{j-1}), \quad 1\leqslant j\leqslant n,\quad \Psi_0:=0.\label{angsys2}
\end{align}

\subsubsection{Unperturbed polar system, part I: radial and angular dynamics}
Since all distances are positive, we consider the part of the sphere with positive coordinates, 
\[
S^+:=S^+(n):=\left\{\vec{x}=(x_1,x_2,\ldots,x_n)^\intercal\in S^{n-1}: x_i>0~\forall i\right\},
\]
on which the following holds.

\begin{proposition}\label{prop:positivity}
$\Sigma_n>0$ on $\overline{ S^+}$. More specifically, $\Sigma_n>\Psi_{j_0}^2/2$, where the index $j_0$ is such that $\Psi_k=0$ for $1\leqslant k<j_0$ and $\Psi_k>0$ for $j_0\leqslant k\leqslant n$. 
\end{proposition}
Note that the index $j_0$ exists since $(\Psi_1,\ldots,\Psi_n)\in S^{n-1}$ does not vanish.
\begin{proof}
Rewrite $\Sigma_n$ as $\Sigma_n=\Psi_{j_0}^3+\sum_{k=j_0+1}^n\Psi_k^2(\Psi_k-\Psi_{k-1})$ and set $f_{k-1}(\Psi_k):=\Psi_k^2(\Psi_k-\Psi_{k-1})$ as a function in $\Psi_k$ with parameter $\Psi_{k-1}$ which has its minimum (on $[0,1]$) at $\frac{2}{3}\Psi_{k-1}$, i.e. $f_{k-1}(\Psi_k)\geqslant f_{k-1}\left(\frac{2}{3}\Psi_{k-1}\right)=-\frac{4}{27}\Psi_{k-1}^3$. Hence,
\begin{equation}\label{ineq1}
\Sigma_n\geqslant\Psi_{j_0}^3-\frac{4}{27}\sum_{k=j_0+1}^n\Psi_{k-1}^3=\Psi_{j_0}^3-\frac{4}{27}\sum_{k=j_0}^{n-1}\Psi_k^3,\qquad \Psi_k=\frac{2}{3}\Psi_{k-1}.
\end{equation}
Since $\Psi_k=\frac{2}{3}\Psi_{k-1}$ one gets  $\Psi_{k}=\left(\frac{2}{3}\right)^{k-j_0}\Psi_{j_0}$ by recursion and consequently,
\begin{align*}
    \Sigma_n&\geqslant\Psi_{j_0}^3-\frac{4}{27}\sum_{k=j_0}^{n-1}\left(\frac{2}{3}\right)^{3(k-j_0)}\Psi_{j_0}^3>\Psi_{j_0}^3\left(1-\frac{4}{27}\sum_{k=0}^{\infty}\left(\frac{2}{3}\right)^k\right)=\frac{15}{27}\Psi_{j_0}^3>\frac{1}{2}\Psi_{j_0}^3
\end{align*}
\end{proof}
\begin{lemma}\label{lem:small}
$\Psi_j'>0$ for $0<\Psi_j\ll 1$ and there exist $T>0, \varepsilon>0$ such that $\Psi_j(t)\in (\varepsilon,1-\varepsilon)$ for all $t\geqslant T$.
\end{lemma}
\begin{proof}
If $\Psi_j=0$, we have $\Psi_{j-1}-\Psi_j+\Sigma_n=\Psi_{j-1}+\Sigma_n>0$. Choosing $\Psi_j>0$ small enough, we therefore get $\Psi_j'=\Psi(\Psi_{j-1}-\Psi_j+\Sigma_n)>0$. This, together with Prop.~\ref{prop:positivity}, implies the second statement by contradiction.
\end{proof}
As a consequence, the distances of the distances $\td_j$ in the unperturbed ODE system are bounded.
\begin{proposition}\label{p:tdbounds}
Solutions to \eqref{e:unperturbed} satisfy 
$\lvert \td_j- \td_k\rvert>0$ and $\td_j-\td_k=\mathcal{O}(1)$ as $t\to\infty$.
\end{proposition}
\begin{proof}
The distances are eventually ordered (cf. Thm.~\ref{thm:orderingunperturbed}), thus $\lvert \td_j-\td_k\rvert>0$ for large $t>0$. To prove the boundedness, suppose that $\td_j-\td_k\to\infty$. Then, $e^{-(\td_j-\td_k)}=\frac{z_j}{z_k}=\frac{\Psi_j}{\Psi_k}\to 0$. However, due to Lemma~\ref{lem:small}, $\frac{\Psi_j}{\Psi_k}\in\left(\frac{\varepsilon}{1-\varepsilon},\frac{1-\varepsilon}{\varepsilon}\right)$ for large $t>0$.
\end{proof}

\subsubsection{Consequences for the perturbed polar system}\label{s:perturbed}
We now return to the system \eqref{radialperturbeddyn}-\eqref{angularperturbeddyn2}. On the sphere (i.e. for $r=0$), the angular equations \eqref{angularperturbeddyn1} and \eqref{angularperturbeddyn2} coincide with those of the unperturbed polar system \eqref{angsys2}. Likewise, near the sphere (i.e. for $0<r\ll 1$), the radial dynamics \eqref{radialperturbeddyn} are dominated by the radial term $-\Sigma_n r$, cf.\ \eqref{rdynamics}. This allows for a perturbation argument from which one can deduce that the distances in the perturbed system diverge. 

\begin{theorem}\label{thm:distances}
There is $\td_0>0$ such for all solutions to \eqref{perturbed} with $\td_{\min}(0)>\td_0$, the statements of Theorem~\ref{thm:orderingunperturbed} and Proposition~\ref{p:tdbounds} hold true.
In particular, there exists $\td_1\geq 0$ such that for all solutions to \eqref{perturbed} with $\td_{\min}(0)>\td_1$ it holds that  $\td_{\min}(t)\geq \td_{\min}(0)$ for all $t\geq 0$.
\end{theorem}

Note that together with Proposition~\ref{prop:Ei2} this in particular implies Theorem~\ref{t:valid}.

\begin{proof}
With respect to the unperturbed polar system \eqref{rdynamics}-\eqref{angsys2}, $\calM_0:=\{r=0\}\times S^+$ is an inflowing normally hyperbolic invariant manifold (cf.\ \cite{HPS,Wiggins}) since the transversal eigenvalues are strictly negative, cf.\ Prop.~\ref{prop:positivity}, and the boundary is repulsive due to Lemma~\ref{lem:small}. By robustness of such invariant manifolds, $\calM_0$ possesses a non-trivial local basin of attraction $\Gamma$ for the unperturbed \emph{and} perturbed polar systems \eqref{radialperturbeddyn}-\eqref{angularperturbeddyn2}.  Theorem~\ref{thm:orderingunperturbed} implies that any solution to the unperturbed system with angular initial data in $S^+$ enters $\Gamma$ in finite time, which therefore also holds for the perturbed system if initial distances are sufficiently large. These solutions thus converge to $\calM_0$. The property that eventually $\Psi_j\in (\varepsilon,1-\varepsilon)$, cf.\ Lemma~\ref{lem:small}, is likewise structurally stable and thus remains valid for the perturbed polar system. In particular, this property means that all distances $\td_j$ diverge. 

Without loss of generality, we can choose $\Gamma$ to lie in the $n$-dimensional local stable manifold $W^s$ of $\calM_0$ in the unperturbed system, which we write as $W^s(S^+)$ since $\calM_0$ is trivially parameterized by $S^+$. Now, for the perturbed and unperturbed polar system, $W^s(S^+)$ is foliated by one-dimensional strong stable fibers $W^{ss}(\Psi)$, respectively, 
\begin{equation}\label{foliation}
W^s(S^+)=\bigcup_{\Psi\in S^+}W^{ss}(\Psi),
\end{equation}
which are pairwise disjoint and each intersect $\calM_0$ in their base point $\Psi\in S^+$. The dynamics of the base points is given by \eqref{angsys2} for \emph{both} the perturbed and unperturbed polar systems, but the fibers differ in general.

The key point is that the flows of both \eqref{rdynamics}-\eqref{angsys2} and \eqref{radialperturbeddyn}-\eqref{angularperturbeddyn2} in $\Gamma$ are slaved to the base point flow such that the perturbed flow inherits the properties of the unperturbed flow as claimed. More specifically, let $\Psi(t)$, $\Psi(0)\in S^+$ be a solution to \eqref{angsys2}. Then the perturbed and unperturbed flows map their respective fibre $W^{ss}(\Psi(0))$ into their respective fibre $W^{ss}(\Psi(t))$; using the foliation \eqref{foliation} for local coordinates of the perturbed system near the sphere, the base flow (which is always the same) decouples from the transverse fibre flow, which differs between perturbed and unperturbed flow. 

Since the base point flow leads to an eventually ordering of the distances for the unperturbed polar system, cf.\  Thm.~\ref{thm:orderingunperturbed}, the slaving implies the same for the perturbed system \eqref{perturbed}. The remaining claims follow analogously. These imply that there is $\td_*$ such that $\td_{\min}(t)$ grows strictly if $\td_{\min}\geq \td_*$, which implies the existence of an $\td_1\leq \td*$ as claimed.
\end{proof}

\subsubsection{Unperturbed polar system, part II: local stability}
Here we add more details to the dynamics of the unperturbed polar system: we prove that there is a unique equilibrium point in $S^+$, and it is locally exponentially stable. For illustration, let us first consider the simplest cases $n=2$ and $n=3$. 

For $n=2$, the system \eqref{angsys2}, with $\Sigma_2=\Psi_1^3-\Psi_1\Psi_2^2+\Psi_2^3$, reads
\begin{align*}
    \Psi_1'&=\Psi_1(-\Psi_1+\Sigma_2)\\
    \Psi_2'&=\Psi_2(\Psi_1-\Psi_2+\Sigma_2).
\end{align*}
Let us first assume that $\Psi_1\neq 0$ which implies $\Psi_1=\Sigma_2$ for an equilibrium. If $\Psi_2\neq 0$, then $\Psi_2=2\Psi_1$ and thus $\Psi_1=\pm\sqrt{\frac{1}{5}}$, and yields the two equilibria $E_{1,2}:=\pm \left(\sqrt{\frac{1}{5}},2\sqrt{\frac{1}{5}}\right)$. For $\Phi_2=0$ we find the two equilibria $E_{3,4}:=(\pm 1,0)$ and for $\Phi_1=0$ the last two, $E_{5,6}:=(0,\pm 1)$. In total, the system has $S_2:=2(2^2-1)=6$ equilibria.

For $n=3$, the situation is already a bit more complicated. We have $\Sigma_3=\Psi_1^3-\Psi_1\Psi_2^2+\Psi_2^3-\Psi_2\Psi_3^3+\Psi_3^3$ and need to distinguish the following different cases depending on the number $0\leqslant k\leqslant 2$ of zero coordinates, each giving a pair of equilibria:
\begin{enumerate}
    \item[(i)] $k=0: E_{1,2}=\pm\left(\sqrt{\frac{1}{14}},2\sqrt{\frac{1}{14}},3\sqrt{\frac{1}{14}}\right)$
\item[(ii)] $k=1: E_{3,4}=\pm\left(\sqrt{\frac{1}{5}},2\sqrt{\frac{1}{5}},0\right),E_{5,6}=\pm\left(\sqrt{\frac{1}{2}},0,\sqrt{\frac{1}{2}}\right), E_{7,8}=\pm\left(0,\sqrt{\frac{1}{5}},2\sqrt{\frac{1}{5}}\right)$ 
\item[(iii)] $k=2: E_{9,10}=(\pm 1,0,0),E_{11,12}=(0,\pm 1,0),E_{13,14}=(0,0,\pm 1)$
\end{enumerate}
In total, the system has $S_3:=2(2^3-1)=14$ equilibria.

For general $n\in\mathbb{N}$ the following holds.
\begin{proposition}
For $n\in\mathbb{N}$, system \eqref{angsys2} posesses $S_n:=2(2^n-1)$ equilibria. Specifically, there are (i) $2n$ equilibria which have exactly one non-zero component (which is therefore $\pm 1$) and (ii) $S_n-2n-2$ equilibria with $2\leqslant j\leqslant n-1$ non-zero components. 
\end{proposition}
\begin{proof}
Let $n\in\mathbb{N}$ and $S_n$ denote the number of equilibria $\vec{\Psi}=(\Psi_1,\ldots,\Psi_n)$ of dimension $n$. First note that $\Psi_j\neq 0$ implies for equilibria $\Psi_j=\Sigma_n+\Psi_{j-1}$ by \eqref{angsys2}. Thus, $\Psi_{j-1}=0, \Psi_j\neq 0$ implies $\Psi_j=\Sigma_n$, and adjacent non-zero entries come as a sequence $(\Sigma_n,2\Sigma_n,\ldots,k\Sigma_n)$, where $2\leq k\leq n$. Hence, an intermediate zero coordinate between two non-zero coordinates results in a triple $(\Sigma_n,0,\Sigma_n)$. 

Consequently, if $\Psi_1=0$, the $S_{n-1}$ equilibria of dimension $n-1$ occur with shifted index; this in particular includes the vectors with $\Psi_j=\pm 1$ for some $2\leq j\leq n$ and zero coordinates otherwise. 
Denote the number of remaining equilibria with $\Psi_1\neq 0$ by $M_n$. By the discussion above, $M_n$ is the number of possibilities to have $k\in\{0,1,\ldots,n-1\}$ zero entries on the positions different from $\Psi_1\neq 0$ (in particular including $(\pm 1,0,0,\ldots,0)$). This number is given by $M_n=2\sum_{k=0}^{n-1}\binom{n-1}{k}$.

All together, 
\begin{equation*}
S_n=M_n+S_{n-1}=\sum_{j=1}^n M_j=2\sum_{j=1}^n\sum_{k=0}^{n-1}\binom{j-1}{k}=2\sum_{j=1}^n2^{j-1}=2(2^n-1).
\end{equation*}
\end{proof}

By the previous lemma,  $E_\pm:=(\Psi_1,2\Psi_1,3\Psi_1,\ldots,n\Psi_1)$ with $\Psi_1=\pm\sqrt{\frac{6}{n(n+1)(2n+1)}}$ are the unique equilibria with non-zero coordinates only. 

\begin{figure}[H]
    \centering
    \begin{tabular}{cc}
     \includegraphics[scale=0.35]{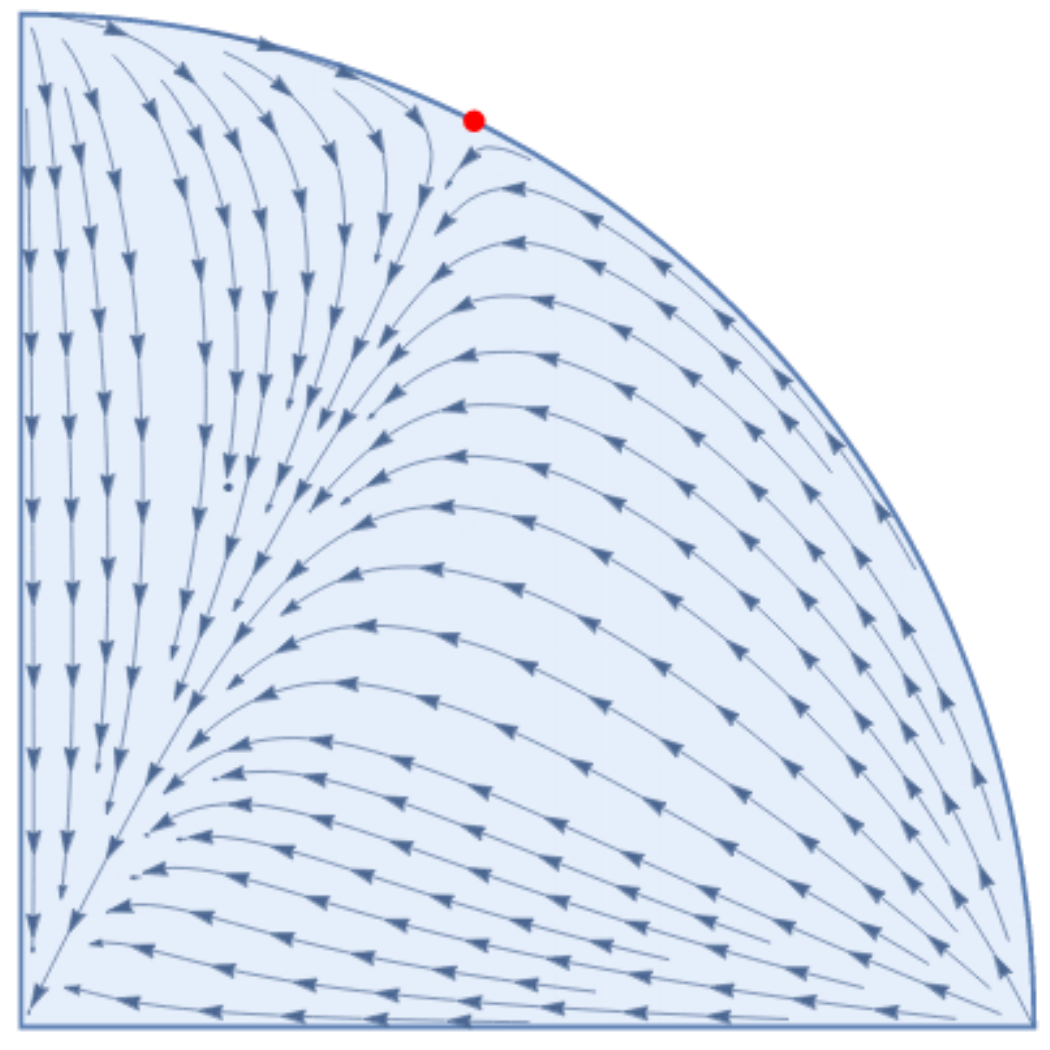} & \includegraphics[scale=0.35]{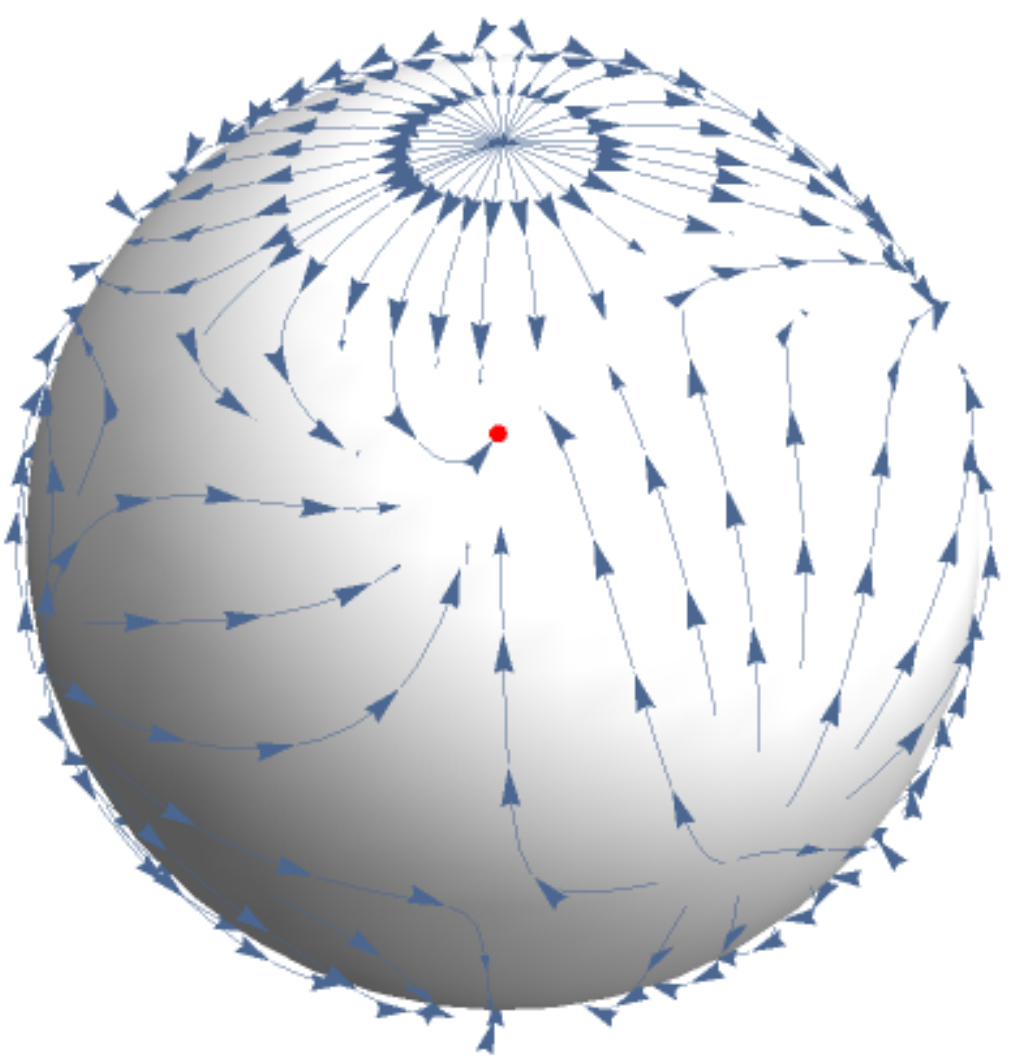}\\
     (a) & (b)
    \end{tabular}
    \caption{Phase portraits and equilibrium $E$ (marked in red) for (a) $n=2$ and (b) $n=3$, respectively. For $n=2$, the phase portrait is restricted to the upper right 
    sector and the transverse direction is neglected since it is not relevant for the angular dynamics. 
    As these simulations suggest, $E$ is a global attractor on $S^+$ for these cases.}
    \label{fig:spheres}
\end{figure}

\textbf{Local stability of $E_+$.} 
In the following, we focus on $E:=E_+$ since for the question of divergence of the distances in the perturbed ODE-system this is the only relevant asymptotic state. In fact, we expect that it is a global attractor on $S^+$ as can be illustrated by simulations for $n=2$ and $n=3$, cf. Fig.~\ref{fig:spheres}. However, it seems difficult to prove this rigorously in general. Instead, by linearizing \eqref{angsys2} in $E$ and determining the eigenvalues of the Jacobian, we show that $E$ is locally stable on $S^{n-1}$.

\begin{theorem}\label{thm:localstability}
The equilibrium $E$ is locally exponentially stable on $S^{n-1}$. The eigenvalues of the linearisation in $E$ are given by $-\frac{k}{6}n(n+1)(2n+1), 2\leqslant k\leqslant n$. For the full unperturbed polar system in $\mathbb{R}^n$, $E$ possesses the unstable eigendirection transverse to $\calM_0$ spanned by $(1,2,\ldots,n)^\intercal$ with corresponding eigenvalue $\frac{1}{3}n(n+1)(2n+1)$.
\end{theorem}

\begin{proof}
The somewhat lengthy proof 
is given in Appendix~\ref{appendixB}.
\end{proof}

\section{Bounded initial kink-antikink data and their annihilation}\label{boundedkinkantikink}

\begin{figure}[h]
    \centering
    \includegraphics{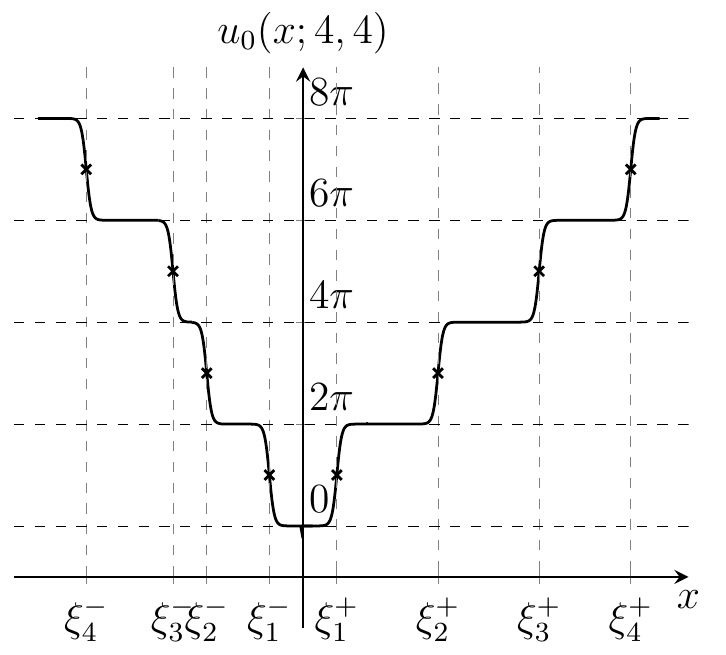}
    \caption{Sketch of initial datum \eqref{data} with $m=n=4$.}
\label{fig:doublestaircase}
\end{figure}

Having discussed pure kink or antikink initial data, we now turn to initial data which are composed of kinks and antikinks, that is, initial data \eqref{data} with $n+m<\infty$ and $nm\neq 0$, cf. Fig.~\ref{fig:doublestaircase} for illustration. Here we constrain ourselves to bounded data and consider the unbounded case in \S\ref{sec:unbounded}. We aim to infer information about the process by which the respective inner pair of kink and antikink with positions $\xi_1^\pm(t)$ annihilate each other.

\subsection{Annihilation process}
In case $m=n$, the equal number of pairs of kinks and antikinks completely annihilate each other in the sense that the solution converges to the rest state $2\pi n$ given by the asymptotic state of the initial data.
\begin{proposition}\label{prop:uniformconvergence}
Let $u_0(x,n,n)$  be an initial datum \eqref{data} with $m=n, 0<n<\infty$. Then, $\lim_{t\to\infty}u(\cdot,t)=2\pi n$, where the convergence is uniform in $x\in\mathbb{R}$. In particular, $\omega(u_0(x;n))=\Omega(u_0(x;n))=\{2\pi n\}$. 
\end{proposition}
\begin{proof}
Consider subsolutions with initial data the pure kink- or antikink sequence $u_0(x;n,0)$, $u_0(x;0,n)$ built from the kink- or antikink positions of $u_0(x;n,n)$, respectively. According to Theorem 1.1 \cite{polacik2016propagating} the kinks and antikinks move with uniform positive, respectively negative speed. The comparison principle implies the claim.
\end{proof}

The following proposition describes the corresponding annihilation process in some more detail by showing that the kinks-antikink pairs get annihilated successively from ``bottom to top'' as expected intuitively.
\begin{proposition}\label{ctimes}
Let $u_0(x;n,n)$ be an initial datum \eqref{data} with $n<\infty$ and corresponding solution $u(x,t)$. Then, there are unique times $0<t_1<t_2<\ldots t_n$ such that $\xi_j^-(t_j)=\xi_i^+(t_j)$, $1\leqslant j\leqslant n$,  i.e., the successively innermost kink and antikink collide at these times. Moreover, there are unique times $t_j^A\in (t_j,t_{j+1})$, $1\leqslant j\leqslant n-1$, such that $\mathrm{argmin}_{x\in\R}\{u(x,t_j^A)\}=2\pi j$, i.e., at these times the successively innermost kink and antikink annihilate.
\end{proposition}
To ease notation we set $t_0^A:=0$.
\begin{proof}
As shown in the proof of Proposition~\ref{prop_posfunc}, for all $t>0$, the solution $u(x,t)$ is monotonically decreasing on $(-\infty,\eta(t))$ with $\lim_{x\to-\infty}u(x,t)=2\pi n$ and monotonically increasing on $(\eta(t),\infty)$ with $\lim_{x\to\infty}u(x,t)=2\pi n$. Together with  Proposition~\ref{prop:uniformconvergence}, this shows that the positions $\xi_i^\pm(t), 1\leqslant i\leqslant n$, collide at some unique time $t_i>0$, respectively, i.e., $u(\eta(t_i),t_i)=(2i-1)\pi$ and $u(\cdot, t;n)>(2i-1)\pi$ for $t>t_i$; therefore $\xi_i^\pm(t)$ exist for $[0,t_i]$ only. In particular, $t_i<t_{i+1}$. Analogously, we identify the annihilation times $t_i^A \in (t_i,t_{i+1})$.
\end{proof}

Next we show that during the annihilation process of the innermost kink-antikink pair, the distances between the remaining kinks (antikinks) satisfy a uniform lower bound. 

\begin{proposition}\label{prop:bounds}
Let $u(x,t)$ be the solution corresponding to an initial condition as above. Then, $d_j^\pm(t)\geq d_{\min}^\pm(0)$ and if $j^\pm$ are indices such that $d_{\min}^\pm(0)=d_{j^\pm}(0)$ then $d_{j^\pm}(t)$ are non-decreasing as long as defined, i.e., $t\leq t_{j}$, respectively.
In particular, $\lvert\xi_{j+1}^-(t)-\xi_{j+1}^+(t)\rvert\geqslant d_{\min}^+(0) + d_{\min}^-(0) + \lvert\xi_{j}^-(t)-\xi_{j}^+(t)\rvert$ for all $t\in(0,t_j]$.
\end{proposition}
\begin{proof}
Due to reflection symmetry, it suffices to prove  the claims for the case `+'. In the proof of Proposition~\ref{prop_non-equi}, we have constructed initial conditions $\underline{u}_{0}$ and $\bar{u}_{0}$ whose corresponding solutions $\underline{u}(x,t)$, $\bar{u}(x,t)$ are sub- and supersolutions, which imply the non-decrease of the minimal distance. In the present case the same subsolutions can be used. However, $\bar{u}_0$ are not providing supersolutions on $\R$, but only on $(-\infty,\eta(t))$ since $u(\cdot,t)$ is minimal at $\eta(t)$ and $\bar{u}(x,t)$ can be above $u(x,t)$ only for $x>\eta(t)$. 

The claims now follow from the statement of Proposition~\ref{prop_non-equi}.
\end{proof}

It is tempting to suppose that the minimal distance between kinks or antikinks, respectively, that arises after annihilations can be used as a lower bound. However, it seems difficult to construct super- and subsolutions to substantiate this.

\medskip
Let us now turn to the case of initial data with unequal number of kinks and antikinks, $m\neq n$; without loss of generality $m<n$. In this situation, the first $\min\{m,n\}$ innermost kink-antikink pairs annihilate and a stacked front of either kinks or antikinks remains after some finite time.  

\begin{proposition}\label{prop:annihilation}
Let $u_0(x;m,n)$ be an initial datum \eqref{data} with $m<n<\infty$ and $mn\neq 0$. The kinks and antikinks at $\xi_i^\pm(t)$, $1\leqslant i\leqslant m$, collide and annihilate in the sense of Proposition \ref{ctimes}. Moreover, the limit sets are given by 
\begin{align*}
\omega(u_0(x;m,n))&=\{2\pi n\},\\
    \Omega(u_0(x;m,n))&=\{2\pi j: m\leqslant j\leqslant n\}\cup\{\varphi_j^+(\cdot-\nu): m<j<n, \nu\in\mathbb{R}\}.
\end{align*}
\end{proposition}
\begin{proof}
In the course of the annihilation process of the $m$ innermost kink-antikink pairs, $u(x,t)$ gets uniformly close to $2\pi m$ for $x\leq \eta(t))$. 
Meanwhile, the distances $d_j^+~(m\leq j<n)$ of the remaining kinks are bounded from below by Proposition~\ref{prop:bounds}, which means the solution gets arbitrarily close to a propagating terrace. 
Since $2\pi m$ is a stable steady state, the solution $u(x,t)$ converges to a propagating terrace and \cite[Theorem 1.1]{polacik2016propagating} implies the statement. 
\end{proof}

\begin{figure}[H]
    \centering
    \includegraphics{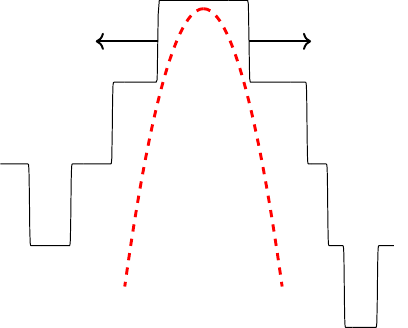}
    \caption{Sketch of initial datum with local maximum and subsolution (red). Kink and antikink which built the local maximum cannot collide.}
\label{fig:multiplemax}
\end{figure}
Towards a more complete picture, let us briefly consider initial data with local maxima built by pairs of kinks and antikinks, cf.\ Fig.~\ref{fig:multiplemax}. For each such maximum with sufficiently large distance of kink and antikink, one can construct a stationary subsolution using phase plane analysis for the equation $u_{xx}+f(u)=0$;
 this shows that kinks and antikinks cannot annihilate at local maxima, but only between local minima as described before. In particular, the limit sets $\omega(u_0)$ and $\Omega(u_0)$ are completely determined by the numbers $m,n$ between maxima.

\section{Unbounded kink or kink-antikink data}\label{sec:unbounded}
\begin{figure}[H]
    \centering
    \begin{tabular}{cc}
   \hspace{-1cm} \includegraphics[scale=0.8]{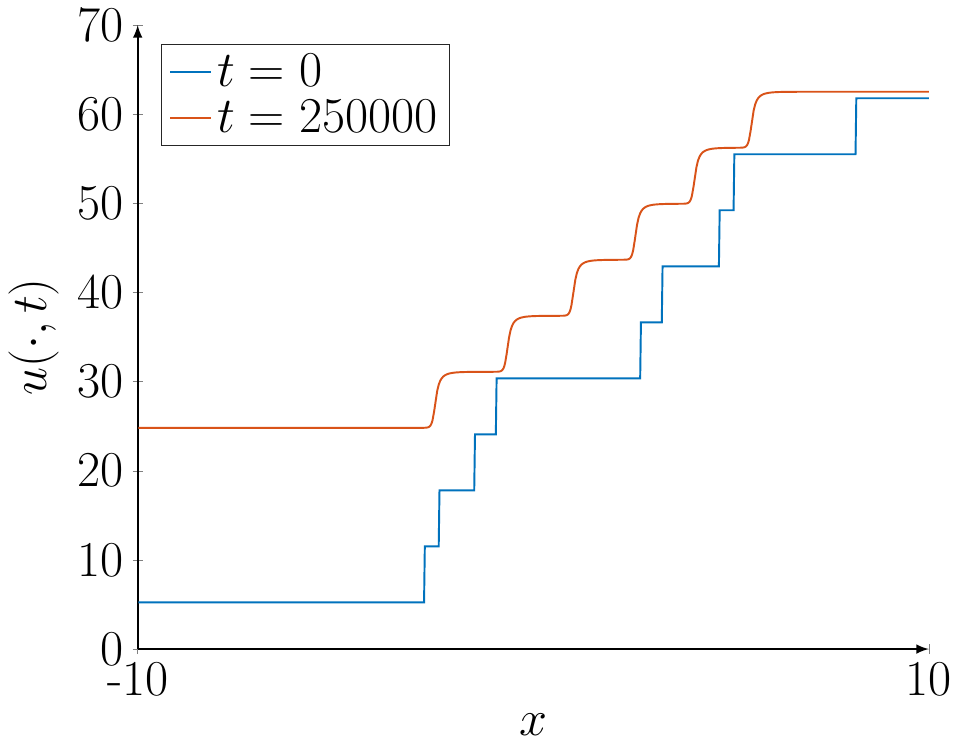} &
\hspace{-0.3cm}\includegraphics[scale=0.8]{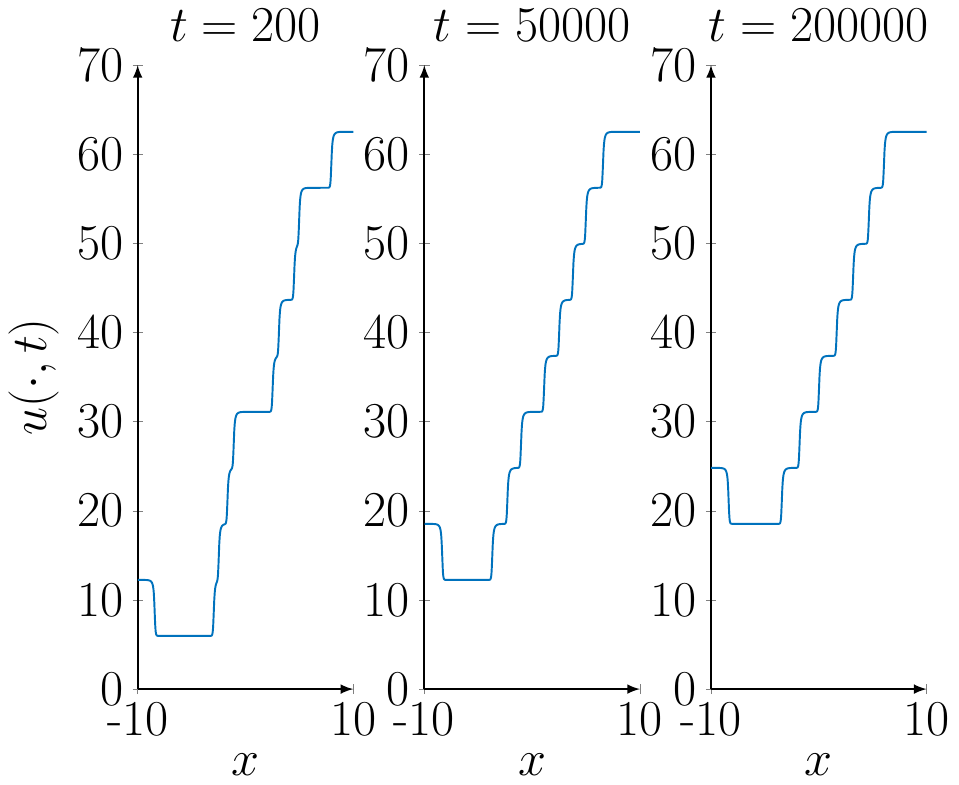}\\
    (a) & (b)
    \end{tabular}
   \caption{Simulations to illustrate the loss of information ($\muf=0.2$). (a) Information encoded in the initial kink distances is locally ``washed out'' at $t=250000$, independent of three `fast' annihilations due to incoming antikinks as plotted in (b). See Appendix~\ref{app:sim} for details on the implementation.} 
    \label{fig:washingout}
\end{figure}
Unbounded initial data is most relevant for the comparison with the dynamics of the cellular automaton GHCA, such as its non-wandering set dynamics and topological entropy. However, for unbounded data several useful results cannot be directly applied, e.g., those on terraces in \cite{polacik2016propagating}. 

Nevertheless, the zero number argument is still applicable (the growth condition of \cite[p. 80]{angenent1988zero} is satisfied for $u\in X_\omega$) so that we readily infer an analogon to Proposition~\ref{ctimes} in case $m=n=\infty$: we can repeatedly choose initial data composed of finite kink-antikink pairs as subsolutions and obtain an infinite sequence of strictly increasing collision and annihilation times. 

However, numerical simulations of \eqref{theta-eq} (lifted to $u(x,t)\in\R$), cf.\ Figs.~\ref{fig:washingout}, \ref{fig:kinkdistPBC}, suggest that  distances equilibrate asymptotically in time.  Hence, initial distance information encoded far from an eventual collision is lost over time, is ``washed out'' -- at least on the scale of the initial data. Indeed, this is consistent with our results on the dynamics of analytic positions for monotone data, and large initial distances. Therefore, as suspected from weak interaction induced by diffusion, we cannot expect entropy and dynamics for \eqref{theta-eq} with unbounded initial data are directly analogous to GHCA.

We next first corroborate the equilibration of distances for unbounded initial data by considering periodic boundary conditions and show that all solutions converge locally uniformly to equidistant staircases. Second, we discuss implications for complexity measures based on positional dynamics.

\subsection{Periodic boundary conditions}

Unbounded superpositions of equidistant kinks (and, analogously, antikinks) are parametrized by the uniform distance $\ell>0$, and the problem is transformed into a boundary value problem under periodic boundary conditions, up to phase rotation. 
\begin{figure}[H]
    \centering
    \begin{tabular}{cc}
    \hspace{-1cm}\includegraphics[scale=0.8 ]{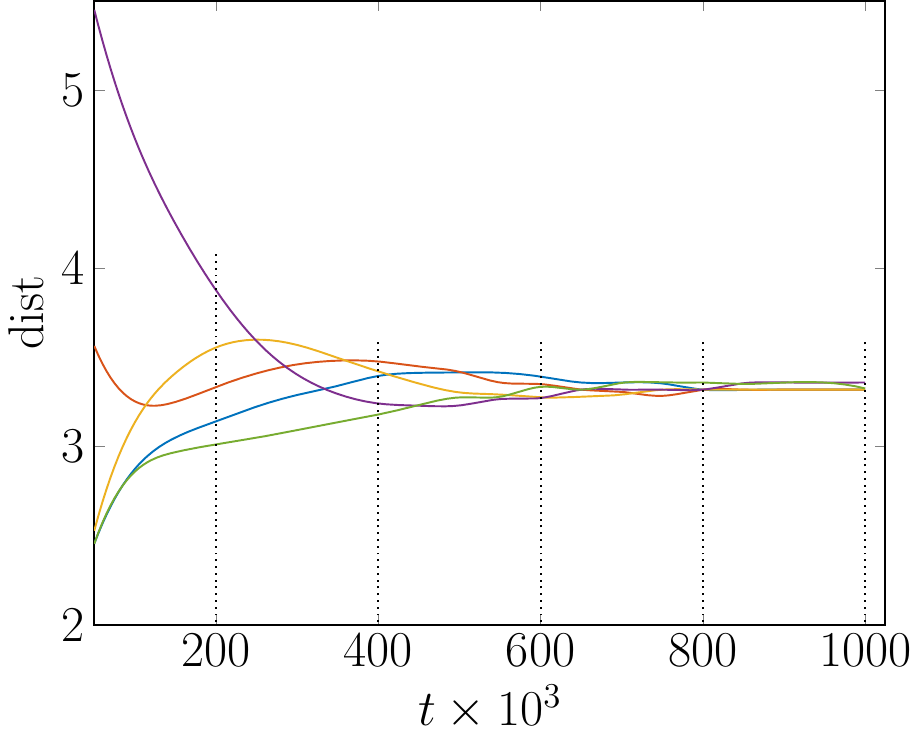} & \hspace{-0.3cm}\includegraphics[scale=0.8]{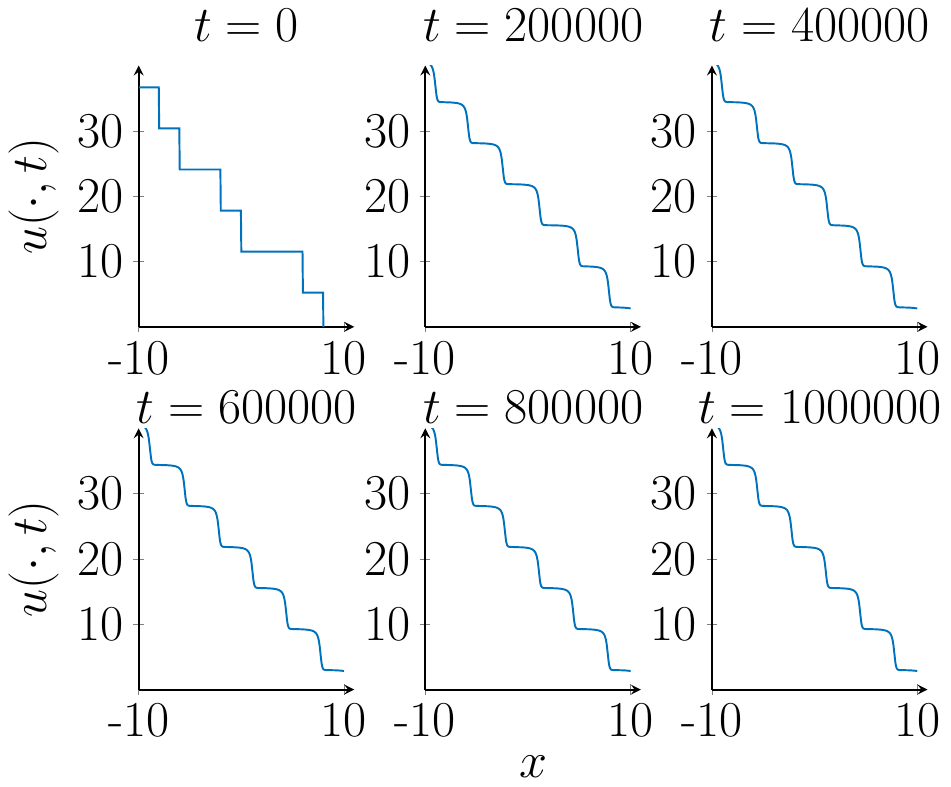}\\
    (a) & (b)
    \end{tabular}
   \caption{Simulation of (a) kink distances and (b) corresponding solutions done with pde2path by freezing under periodic boundary conditions and with $\muf=0.2$, cf.\ Appendix~\ref{app:sim}. The initial datum converges to an equidistant state.}
\label{fig:kinkdistPBC}
\end{figure}

More precisely, for given $\ell\in\R_+$ and $j\in\mathbb{N}$, we consider the initial-boundary value problem 
\begin{subnumcases}{\label{IBVP}}
\theta_t=\theta_{xx}+f(\theta), & $-\ell\leqslant x\leqslant \ell$ \label{IBVP1} \\
\theta(0,x)=\theta_0(x), & $ -\ell\leqslant x\leqslant\ell$ \label{IBVP2}\\
\theta(-\ell)=\theta(\ell)+2\pi j, & $\theta(-\ell)\leqslant \theta\leqslant\theta(\ell)$\label{IBVP3}\\
\theta'(-\ell)=\theta'(\ell)\label{IBVP4}
\end{subnumcases}
and prove in Prop.~\ref{p:washout} below that solutions of \eqref{IBVP} converge to a unique equidistant staircase in terms of the $\omega$-limit set of $\theta_0$.

To this end, we first consider travelling wave solutions $\theta(x,t)=u(x-at)=:u(z)$ and focus on the following boundary value problem, where $'=\frac{d}{dz}$.
\begin{proposition}\label{prop:bvp}
For each triple $(\ell,j,a)\in \mathbb{R}^+\times\mathbb{N}\times\mathbb{R}$, the boundary value problem
\begin{subnumcases}{\label{TWE}}
       u''+a u'+f(u)=0, & $-\ell<x<\ell$ \label{TWE1}\\
        u(-\ell)=2\pi j, u(\ell)=0, & $0\leqslant u\leqslant 2\pi j$\label{TWE2}
    \end{subnumcases}
has a unique, monotonously decreasing solution $u$ with $\partial_x u<0$ on $(-\ell,\ell)$. Moreover, there exists a unique $a(\ell,j)$ such that the unique solution $u$ of \eqref{TWE} corresponding to $(\ell,j,a(\ell,j))$ additionally satisfies $u'(-\ell)=u'(\ell)$. This $a(\ell,j)$ is given by 
\begin{equation*}
a(\ell,j)=(F(2\pi j)-F(0))\left(\int_{-\ell}^\ell (u')^2\, dx\right)^{-1},
\end{equation*}
where $F'=f$. In particular, $\lim_{\ell\to 0}a(\ell,j)=0$.
\end{proposition}

We remark that the case $j=1$ is essentially contained in \cite[Thm 3.5]{ErmentroutRinzel}.

\begin{proof}
Sub- and super solutions are given by $\underline{u}\equiv 0$ and $\bar{u}\equiv 2\pi j$, respectively, which shows the existence of a solution $u$. Applying Theorem 1.4 \cite{Berestycki1991} (``sliding method'') to $\Omega_\ell:=[-\ell,\ell]$, together with the maximum principle, implies the monotonicity and uniqueness of this solution. 

Moreover, this unique solution (i) depends continuously on $a$ and (ii) is strictly decreasing in its dependence on $a$ (i.e. if $a_1<a_2$, then $u_2<u_1$ in $\Omega_{\ell}$, where $u_1$ and $u_2$ are the solutions of \eqref{TWE} corresponding to $(\ell,j,a_i), i=1,2$, cf. Corollary 5.1 \cite{Beresticky1992497}). 
By Lemma 5.2 \cite{Beresticky1992497}, one infers that (iii) $\lim_{a\to -\infty}u=2\pi j$ and $\lim_{a\to\infty}u=0$, both uniformly in $x$.  
Items (i)-(iii) together imply the existence of a unique solution $u$ and a unique $a(\ell,j)$ such that $u'(-\ell)=u'(\ell)$.

As to $a(\ell,j)$ and its asymptotics, one multiplies equation \eqref{TWE1} by $u'$ and integrates to get   
\begin{equation*}
    a(\ell,j)\int_{-\ell}^\ell (u')^2\, dx=-\int_{-\ell}^\ell u'u''+u'f(u)\, dx=F(2\pi j)-F(0),
\end{equation*}
where $u'(-\ell)=u'(\ell)$ is used. We obtain $a(\ell,j)=(F(2\pi j)-F(0))\left(\int_{-\ell}^{\ell}(u')^2\, dx\right)^{-1}$. 
By the Cauchy-Schwarz inequality,
\begin{equation*}
    \int_{-\ell}^\ell (u')^2\, dx\geqslant\frac{4\pi^2 j^2}{2\ell}
\end{equation*}
and thus $a(\ell,j)\to 0$ as $\ell\to 0$.
\end{proof}
\begin{remark}
Since $f$ is $2\pi$-periodic, the statement clearly remains true for boundary conditions $u(\ell)=2\pi k$ and $u(-\ell) = 2\pi(k+j), j,k\in\mathbb{N}$,
and the solution is just $u+2\pi k$; in particular  \eqref{IBVP3}, \eqref{IBVP4} hold. We also remark that $F(2\pi j)\neq F(0)$ for our nonlinearity $f$.
\end{remark}

In the following, we focus on the unique solution $u$ of \eqref{TWE} corresponding to the triple $(\ell,j,a(\ell,j))$. Due to its additional property $u'(-\ell)=u'(\ell)$, it is the relevant solution for the purposes of this section. As a consequence of the previous proposition, its shape is determined by single periodic solutions in the following sense; in particular all kinks in $u$ are equidistant.
\begin{proposition}\label{prop:bvp1}
Let $u$ denote the unique solution of \eqref{TWE} which corresponds to the triple $(\ell,j,a(\ell,j))$. Then $u$ consists of $j$ space shifted and phase rotated copies of the solution $\tilde{u}$ of \eqref{TWE} corresponding to $(\tilde{\ell},1,a(\tilde{\ell},1))$ with $\tilde{\ell}=\ell/j$. In particular, $a(\ell,j)=a(\tilde{\ell},1), u'(\pm\ell)=\tilde{u}'(\pm\tilde{\ell})<0$ and $u$ has time period $T=\tilde{\ell}/a(\tilde{\ell},1)$.
\end{proposition}

\begin{proof}
We split the interval $[-\ell,\ell]$ into $j$ intervals $I_s:=[\ell_{s+1},\ell_s],~(0\leqslant s\leqslant j-1)$ of width $2\tilde{\ell}$, where $\ell_s:=\ell-2s\tilde{\ell}$. In particular, $\ell_{s}>\ell_{s+1}$ and $\ell_0=\ell, \ell_j=-\ell$. On each interval, we consider a space shifted and phase rotated version of $\tilde{u}$ which together build a solution due to equal derivatives at the boundaries. More precisely, let $v(x):=\tilde{u}(x-(j-1)\tilde{\ell})$ be the spatial shift of $\tilde{u}$ to the rightmost interval $I_0$ and, for $s\in\{0,1\ldots j-1\}$, set
\begin{equation*}
    U(x):=v(x+2s\tilde{\ell})+2\pi s,\qquad x\in I_s.
\end{equation*}
The function $U$ is thus defined on $[-\ell,\ell]$ and solves \eqref{TWE} with $a=a(\tilde{\ell},1)$ and, moreover, $U'(-\ell)=U'(\ell)$. 
By the uniqueness result of Prop.~\ref{prop:bvp}, it follows that $U=u$ and, in particular, $a(\ell,j)=a(\tilde{\ell},1)$ and $u'(\pm\ell)=\tilde{u}'(\pm\tilde{\ell})$. By the construction of $U$, the same proposition implies $\tilde{u}'(\pm\tilde{\ell})<0$ since $u'(x)<0$ for all $x\in (-\ell,\ell)$.  
\end{proof}

Having established the unique solution of problem \eqref{TWE} with $u'(-\ell) = u'(\ell)$, 
we next show that solutions of \eqref{IBVP} converge to this solutions in the following sense; without loss of generality, we choose $\theta(\ell)=0$ and consider continuous initial conditions $0\leqslant \theta_0\leqslant 2\pi j, j\in\mathbb{N}$.

\begin{proposition}\label{p:washout}
Let $u$ be a solution of \eqref{IBVP} with $\theta(\ell)=0$ and initial datum $0\leqslant \theta_0\leqslant 2\pi j$ for some $j\in\mathbb{N}$. Then, the $\omega$-limit set of $\theta_0$ consists of the orbit associated to the unique solution of \eqref{TWE} corresponding to $(\ell,j,a(\ell,j))$.
\end{proposition}
\begin{proof}
This is basically a consequence of \cite[Theorem 1]{FiedlerMallet} which characterises the limit set by providing a dichotomy between stationary and periodic solutions. However, in order to apply this theorem, we need to transform \eqref{IBVP} into a problem with periodic boundary conditions.

To this end, we consider  $w(x,t):=\theta(x,t)+\frac{\pi j(x-\ell)}{\ell}$ which transforms the problem into
\begin{subnumcases}{\label{IBVP'}}
       w_t=w_{xx}+g(x,w), \label{IBVP'1}\\
       w(0,x)=w_0(x)=\theta_0(x)+\frac{\pi j (x-\ell)}{\ell}, \label{IBVP'2}\\
       w(-\ell)=w(-\ell), \label{IBVP'3}\\
       w'(-\ell)=w'(\ell) \label{IBVP'4}
\end{subnumcases}
where $g(x,w):=f\left(w-\frac{\pi j(x-\ell)}{\ell}\right)$. In particular, there is a one-to-one relation between solutions of \eqref{IBVP} and \eqref{IBVP'}. Let $\hat{u}$ denote the unique solution of \eqref{TWE} associated with $(\ell,j,a(\ell,j))$. Since $W(x,t)=\hat{u}(x-a(\ell)t)+\frac{\pi j(x-\ell)}{\ell}$ is a periodic solution solution of \eqref{IBVP'}, the mentioned dichotomy implies that $\omega(w_0)$ consists of the periodic orbit associated with $W$ and transforming back proves the statement.
\end{proof}

\begin{remark}The previous proposition shows that solutions of \eqref{IBVP} converge locally uniformly to (some translate of) the corresponding unique solution of \eqref{TWE}, i.e., the kinks become equidistant, cf.\ Fig.~\ref{fig:kinkdistPBC}. For a general analysis of convergence in one-dimensional semilinear heat equations, including other types of boundary conditions and nonlinearites, we refer to \cite{CHEN1989160,Matano,FiedlerAngenent} and references therein. 
\end{remark}

\subsection{Complexity considerations}
The motivation to consider kink-antikink dynamics in the $\theta$-equations emerges from our analysis of GHCA and their (topological) complexity \cite{kessebhmer2019dynamics}. The understanding of this complexity mainly relies on the observation that for GHCA the original and somehow abstract definitions of topological entropy \cite{Adler,bowen1971entropy} can be substantiated by considering simpler but equivalent definitions. In this regard, one combinatorial approach is to count possible realizations of bounded space-time windows and to determine the exponential growth rate as the size of these windows increases \cite{DurSte91,kessebhmer2019dynamics}. On the other hand, the decomposition of the non-wandering set reveals that the topological complexity is completely determined by the invariant subsystem consisting of bi-infinite configurations which are composed of counter propagating pulses,
\begin{equation}\label{sequences}
(\ldots, 0_{k_j}, \P^\Rrm, 0_{k_{j+1}},\ldots,\P^\Rrm,0_{k_0},\P^\Lrm,0_{k_1},\P^\Lrm,\ldots),
\end{equation}
where $0_k:(0,\ldots,0)$ are zero blocks of length $k\in\mathbb{N}$ and $\P^{\Rrm,\Lrm}$ represent local pulses, i.e., blocks that move to the right and left under the CA-dynamics, respectively. Thus, a more tailored way to encode the topological complexity of GHCA relies on the construction of a topological conjugacy by defining a proper homeomorphism that maps counter propagating configurations to admissible sequences of \emph{collision sequences} in the form of pairs $(p_n,s_n)_{n\in\mathbb{N}}\subset\mathbb{Z}\times\mathbb{N}$ of collision positions $p_n$ and times $s_n$. 

We stress that the consistency of these two approaches to determine the topological entropy is crucially based on the specific dynamics and topological setup of GHCA and is far from general validity. 

Nevertheless, for $\theta$-equations, tracking the geometric or analytic kink and antikink positions provides a mapping from kink-antikink initial data to collision sequences and thus allows for a comparison with those of GHCA with its encoded complexity. 
Hence, insight into position dynamics of kink, antikink and kink-antikink initial data is a preqrequisite for at least a first heuristic insight into the underlying complexity of the dynamics

As we have shown, bounded initial data give finite sequences of this kind, but for unequal numbers of kinks and antikinks the positional dynamics remains non-trivial also after the final collision in terms of terraces. This eventually occurs on an exponentially slow time scale and distances of terraces eventually diverge, cf.\ Prop.~\ref{prop_equi}-\ref{prop_non-equi}, so that we do not expect a significant contribution to any position-based complexity measure. This may be corroborated further based on our lower bounds for the distances and the ODE for large initial distances, which severely constrains the positions and thus a priori reduces a position based complexity. This metastable slow dynamics already occurs before collisions and thus modifies the relation of initial positions and collision sequence. However, this is a perturbation for any fixed number of $n$ initial kinks and antikinks, and we conjecture this can be compensated by a perturbation of initial positions.

In contrast, the above result on equilibration of distances strongly suggest that this is no longer true for unbounded kink-antikink initial data. 
While the dynamics of (semi-)unbounded kink or antikink data resembles the pure shift dynamics of (semi-)infinite pulse configurations of GHCA, the dynamics of unbounded kink-antikink data bears resemblance to that of infinite counter propagating pulses of GHCA -- up to the aforementioned equilibriation of distances for `late' collisions.

For a direct comparison to the complexity of GHCA in terms of positions, a more specific question is whether any given admissible sequence of collision positions and times for GHCA can also be realized by appropriate initial kink-antikink data. To be more concrete, let $x=(x_n)_{n\in\mathbb{Z}}$ denote a configuration of the form \eqref{sequences} with pulse positions $p_i^-$ (right-moving) and $p_i^+$ (left-moving), respectively, which realizes a given sequence $(p_n,s_n)_{n\in\mathbb{N}}$ of collision points and times. On the one hand, if we choose kink and antikink positions to be exactly the same, $\xi_i^\pm=p_i^\pm$, in general one cannot expect to observe the given collisions positions and times (cf. ``washing out'', Fig.~\ref{fig:washingout}). On the other hand, this does not rule out the possibility that a time-rescaled sequence $(p_n,\tau s_n)_{n\in\mathbb{N}}$, with suitable scaling factor $\tau$, can be attained from modified initial kink and antikink positions, especially because the equilibration of distances is a phenomenon under large distances, i.e., far from collisions only.  This shows that more general quantitative knowledge about the position dynamics, including less restrictive initial data, is required in order to analyze these aspects rigorously. 

\section{Discussion}\label{sec:disc}
 
Motivated by studies of the Greenberg-Hastings cellular automaton (GHCA) as a caricature of excitable systems, in this paper we have considered the $\theta$-equations describing oscillatory phase dynamics, as the perhaps simplest PDE model of excitable media. Since the non-wandering set of GHCA in essence consists of certain excitation pulse sequences \cite{kessebhmer2019dynamics}, we have focussed on the analogue of such data in $\theta$-equations, which consists of kinks and antikinks. Moreover, in GHCA the topological entropy can be related to kink-antikink collisions. We have therefore analyzed the dynamics of bounded and unbounded kink-antikink initial data, including pure kink and antikink data. To this end, we have defined geometric and analytic positions of kinks and antikinks, the first used for a qualitative analysis of bounded and unbounded data, the latter for quantitative results concerning bounded monotone data. For bounded initial data, the theory of terraces shows that, up to spatial reflection, the $\omega$-limit set indeed consists essentially of finite kink-sequences that weakly interact \cite{polacik2016propagating,Risler}.

 As to the qualitative analysis, we have shown that the set of geometric positions is well-defined 
 and consists of isolated points that lie on smooth curves up to collisions. Using the comparison principle, we have revealed that the minimal initial distance is a global lower bound for distances and that collision times and positions can be tracked abstractly. As a model for unbounded data far from collisions, we have considered monotone data with periodic boundary conditions and have shown that the initial distances asymptotically equilibrate. Consequently, information encoded in these distances is lost over time which indicates that for the PDE, contrary to GHCA, an analogous topological entropy based on positions alone cannot be expected. 
 
 As a quantitative analysis, for bounded monotone data, we have derived ODE for the analytic positions of the kinks (antikinks) within the weak interaction regime. By blow-up type singular rescaling and a perturbation argument, we have shown that the dynamics is slaved to spherical dynamics and distances become ordered in finite time, and eventually diverge; again, this incorporates a loss of information in terms of initial distances. 
 
 The combination of comparison principle and weak interaction theory has revealed that the kink-antikink collision dynamics in the PDE is a multi-scale problem, in contrast to the single scale nature of the GHCA. The fast time scale is essentially determined by the speed of an individual kink, while the slow time scale stems from `tail' interaction that is exponentially slow in the kink (or antikink) distances.
 
In order to combine these approaches for unbounded data, it would be interesting to study whether the approach of \cite{ZelikMielke} can be used to admit infinitely many kinks (antikinks) in order to justify ODE for the positions, and whether the approach of \cite{WrightScheel, Selle} -- which requires different speeds of the kinks (antikinks) -- can be adapted to derive motion laws for non-monotone solutions from kink-antikink initial data.
This might allow to estimate complexity measures adapted to the PDE context, and thus quantify the impact of the slow weak interaction on the fast collision dynamics.

Concerning models of excitable media from systems of PDE, such as the famous FitzHugh-Nagumo equations, there are two major issues. First, the interaction of excitation pulses no longer needs to be pure annihilation, but rebounds and even pulse replication is possible \cite{carter2020pulse,hayase2000self,nishiura2000self}. Second, while weak interaction theory can be generalised to systems, the comparison principle cannot. Nevertheless, the heuristic conclusions extend to this case: weak and strong interaction yield a multi-scale problem and diffusion effects positional dynamics in a non-trivial way for infinitely many pulses, thus impacting positional complexity.

Lastly, we mention that the comparison principle has been replaced by energy methods in the Allen-Cahn and Cahn-Hilliard equations
\cite{Scholtes2018,West2020}. However, typical PDE systems of excitable media do not seem to possess an energy structure that could be exploited.

\bigskip\textbf{Acknowledgements: }\\
This work has been supported by grant RA 2788/1-1 of the German research fund (DFG).

\appendix
\section{Law of motion - ODE for the kink distances}\label{derivationODE}

This appendix is devoted to the derivation of the differential system (\ref{eq:20}). 
As explained in the introduction of \S\ref{s:monotone_ODE}, we follow the scheme 
presented in \cite{Ei2002,Rossides2015} that we adapt to our situation for the convenience of the 
reader.

\paragraph{Preliminaries}
We consider the parabolic problem (\ref{theta-eqz}), \textit{i.e.} already in the moving framework 
$z=x-ct.$ It admits a family of standing wave profiles $\{\vp(\cdot-\xi)+2k\pi,\ \xi\in\R,k\in\Z\},$ 
where $\vp$ is the unique standing wave of (\ref{theta-eqz}) connecting $2\pi=\vp(-\infty)$ 
and $0=\vp(+\infty)$ such that $\vp(0)=\pi.$ Analyzing \eqref{e:kink}, it is a well known fact that $\vp$ satisfies
\begin{equation}\label{eq:41}
\vp(z)=\left\{
\begin{array}{lc}
2\pi-a_-e^{\mu z}\lp1+\mathcal{O}(e^{\gamma z})\rp & z< 0, \\
a_+e^{\lambda z}\lp1+\mathcal{O}(e^{-\gamma z})\rp & z\geq 0,
\end{array}
\right.
\end{equation}
for positive $a_-,a_+,\gamma,$ and where $\lambda,\mu$ are the eigenvalues 
of the linearized system at the asymptotic states $(\vp,\vp')=(0,0)$ and $(\vp,\vp')=(2\pi,0),$ 
cf.\ \eqref{eq:42}.

Linearizing (\ref{theta-eqz}) at these equilibria gives the same linear operator by periodicity of $f$ given by
\begin{equation}\label{eq:43}
\ML : \left\{\begin{array}{ccl}
\mathcal{D}(\ML)\subset\MX & \longrightarrow & \MX \\ 
v & \longmapsto & \dr_{zz}v+c\dr_zv+f'(\vp)v.
\end{array}\right.
\end{equation}
For the purpose of our analysis, one can take $\MX=L^2(\R)$ and $\mathcal{D}(\ML)=H^2(\R).$
Due to the translation of the steady state equation, the operator $\ML$ has a kernel: 
$\ML\vp'=0.$ Since $\vp'$ has constant sign and considering the behavior of $f'(\vp)$ as 
$z\to\pm\infty$ it follows from Sturm-Liouville theory that 0 is the largest eigenvalue of $\ML,$ 
it is simple, and isolated: $\MN(\ML)=\MN(\ML^2)=\vp'\R.$
It is a common result (see \cite{SATTINGER1976312} or \cite{JMR_ARMA2} for instance)
that there exists a closed subspace $\MX_1\simeq\mathcal{R}(\ML)$ of $\MX$ such that 
$\MX=\MX_1\oplus\MN(\ML).$ The space $\MX_1$ is the kernel of an element of the dual $\MX^*$ 
that we denote $e^*,$ with the normalization $\ds \langle \vp',e^*\rangle=1.$ It is defined 
through the usual inner product 
\begin{equation}\label{eq:44}
\langle\psi,e^*\rangle=\Lambda\int_\R e^{cz}\vp'(z)\psi(z)\rmd z,
\end{equation}
where $\Lambda$ is the normalization constant. It follows that the projection onto $\MN(\ML)$
 is given by $P\psi=\langle\psi,e^*\rangle\vp'.$
 
 In our situtation, we are dealing with multiple fronts defined by their positions 
 $\eta_1>\eta_2>\cdots>\eta_n.$ With a straighforward abuse of notation, we denote
 \begin{equation*}
\vp_j=\vp(\cdot-\eta_j),\qquad \ML_j=\dr_{zz}+c\dr_z+f'(\vp_j),\qquad e^*_j\simeq z\mapsto e^{c(z-\eta_j)}\vp'_j(z).
 \end{equation*}
Notice that if the $\eta_j$ depend on time, so do these objects.

\paragraph{Projection Scheme}
Let us now turn to the proof of Lemma~\ref{lem:decomp} and define
the map 
\begin{equation}\label{eq:45}
\Phi:\left\{
\begin{array}{ccl}
U\subset \R^{n} & \longrightarrow & \R^n \\
(\eta_1,\cdots,\eta_n) & \longmapsto & \lp \langle w,e_1^*\rangle,\cdots ,\langle w,e_n^*\rangle  \rp
\end{array}\right.
\end{equation}
where $w=v-\sum_{j=1}^n\vp_j$. The Jacobian
\begin{equation*}
\frac{\partial \Phi}{\partial\eta}=[\langle \varphi_j',e_i^*\rangle-\delta_{i,j}\langle w,e^{c(z-\eta_i)}\varphi_i''\rangle]_{i,j=1}^n
\end{equation*}
has diagonal entries close to 1 for $0<\eps,\delta^{-1}$ small enough, i.e. $w$ small and distances $|\eta_i-\eta_j|$ large enough, while the off-diagonal entries are close to zero. The implicit function theorem gives the desired result.

\paragraph{Reduced ODE}

We turn out to the derivation of (\ref{eq:20}). Let us emphasize once 
again that we refer to \cite{Ei2002} for further details and proofs of all the underlying estimates. 

To get the reduced ODE equation we use the Ansatz \eqref{ansatz}. It yields
\begin{equation}\label{eq:A1}
\dr_tu=-\sum_{j=1}^{n}\eta_j'\varphi_j'+\dr_tw=L[u],
\end{equation}
where $\ds L[v]:=\dr_{zz}v+c\dr_{z}v+f(v)$ denotes the nonlinear operator associated
with our problem. Let us project \eqref{eq:A1} onto the kernel of $\mathcal{L}_j^*,$ 
for all $j\in\{1,\dots,n\}:$
\begin{equation}\label{eq:A2}
-\eta_j'-\sum_{i\neq j}\eta_i'\langle e_j^*,\varphi_i'\rangle+\langle e_j^*,w_t\rangle=\langle e_j^*,L[u]\rangle.
\end{equation}
On the other hand, since $\ds \frac{\rmd }{\rmd t}\langle e_j^*,w\rangle=0,$ it follows that
$\ds \langle e_j^*,w_t\rangle=\eta_j'\Lambda\int_\R e^{c(z-\eta_j)}w(z)\varphi_j''(z)\rmd z.$ 
The system for the positions $\eta=(\eta_1,\ldots,\eta_n)^\intercal$ is therefore 
given by $\ds A\eta=\beta,$ where $\ds \beta_j=\langle e_j^*,L[u]\rangle$ and 
the matrix $A=(a_{i,j})$ is given by 
\begin{align}
a_{i,i}=&-1+\Lambda\int_\R e^{c(z-\eta_j)}w(z)\varphi_j''(z)\rmd z \nonumber \\
a_{i,j}= &\langle e_j^*,\varphi_i'\rangle\textrm{ for }i\neq j. \nonumber
\end{align}
Notice that if $i\neq j,$ $a_{i,j}=O(e^{-\mu\delta_{\min}})$ and that $a_{i,i}=-1+O(\lV w\rV),$
therefore in the expected asymptotics the matrix $A$ is nearly equal to $-I_n.$ The reduced 
ODE for the positions becomes 
\begin{equation}\label{eq:A3}
-\eta_j'(t)=\langle e_j^*(t),L[u(t)]\rangle+\textrm{hot},
\end{equation}
where the higher order terms $\textrm{hot}$ are of order $\lV w\rV_{H^2(\R)}+e^{-\mu\delta_{\min}}\sum\eta_i'.$ In the following, we explain how 
this reduced ODE \eqref{eq:A3} leads to \eqref{eq:20}. For the sake of simplicity, 
during our computation we will include any other negligible terms into `$\textrm{hot}$', without 
changing its name, keeping in mind the leading order we are interested in. We also focus 
on the (most intricate) situation $1<j<n,$ the two specific cases $j=1$ and $j=n$ being 
similar. Let us first notice that from \eqref{ansatz},
\begin{align}
L[u] = & \sum_{i=1}^{n}\lp \varphi_i''+c\varphi_i'\rp+w_{zz}+cw_z+f\lp\sum_{i=1}^{n}\varphi_i+w\rp \nonumber \\
 = & f\lp\sum_{i=1}^{n}\varphi_i\rp-\sum_{i=1}^{n}f(\varphi_i)+O\lp\lV w\rV_{W^{2,\infty}(\R)}\rp \nonumber
\end{align}
so \eqref{eq:A3} becomes 
\begin{equation}\label{eq:A4}
-\eta_j'(t)=\left\langle e_j^*,f\lp\sum_{i=1}^{n}\varphi_i\rp-\sum_{i=1}^{n}f(\varphi_i)\right\rangle+\textrm{hot}.
\end{equation}
The above bracket is an integral over $\R$ (see \eqref{eq:44}). Let us fix $M$ large enough,
such that $0<M<\delta_{\min},$ to be determined later. We split the integration over three domains
\begin{align}
\R= & (-\infty,\eta_{j+1}+M]\cup [\eta_{j+1}+M,\eta_{j-1}-M]\cup[\eta_{j-1}-M,\infty), \nonumber
\end{align}
so that 
\begin{equation}
 -\eta_j'(t)= I_1^-+ I_2+ I_1^++\textrm{hot}.
\end{equation}

\paragraph{Leading order terms: integration around the front} Let us first focus on $I_2,$ 
since it is the most important term.
Notice that for $z\in [\eta_{j+1}+M,\eta_{j-1}-M]$ one has, from \eqref{eq:41} and classical
invariant manifold theory that 
\begin{equation}\label{eq:A5}
\left\{	\begin{array}{cl}
\varphi_i(z) = a_+e^{\lambda(z-\eta_j)}+\textrm{hot}, & \ \textrm{ for }i>j \\ 
2\pi-\varphi_i(z) = a_-e^{\mu(z-\eta_i)}+\textrm{hot}, & \ \textrm{ for }i<j,
\end{array}\right.
\end{equation} 
and that their derivatives satisfy similar estimates. Then, using the $2\pi-$periodicity of $f,$ and linearizing at $\varphi_j$ and $0,$ one has (see also \eqref{eq:A13} and \eqref{eq:A132}):
\begin{align}
\frac{1}{\Lambda}I_2=& \int_{\eta_{j+1}+M}^{\eta_{j-1}-M}e_j^*(z)\lp f\lp\sum_{i=1}^n\varphi_i(z)\rp-\sum_{i=1}^n f(\varphi_i(z))\rp\rmd z \nonumber\\
= &\int_{\eta_{j+1}+M}^{\eta_{j-1}-M}e_j^*\lp f\lp \varphi_j+\sum_{i<j}(\varphi_i-2\pi)+\sum_{i>j}\varphi_i\rp-f(\varphi_j)-\sum_{i<j}f(\varphi_i-2\pi)-\sum_{i>j}f(\varphi_i)\rp\rmd z\nonumber \\
 = & \int_{\eta_{j+1}+M}^{\eta_{j-1}-M}e_j^*(z)\lp f'(\varphi_j(z))-f'(0)\rp\lp \sum_{i<j}(\varphi_i(z)-2\pi)+\sum_{i>j}\varphi_i(z)\rp\rmd z+\textrm{hot} \label{eq:A6}
\end{align}
where the higher order terms in \eqref{eq:A6} are of order
 $\ds \textrm{hot}= O\lp (\delta_j+\delta_{j-1})e^{-2\mu M}\rp$ as can be seen by \eqref{eq:A5}. We split the integral 
 in \eqref{eq:A6} between the left interacting fronts ($i>j$) and the right interacting 
 fronts ($i<j$). Once again, the leading order terms are determined by 
 the closest fronts, $\varphi_{j+1}$ and $\varphi_{j-1}$ respectively. It follows that 
\begin{align}
\frac{1}{\Lambda}I_2 = &  \int_{\eta_{j+1}+M}^{\eta_{j-1}-M}e_j^*(z)\lp f'(\varphi_j(z))-f'(0)\rp \lp\varphi_{j-1}(z)-2\pi\rp \rmd z \nonumber \\ 
&+ \int_{\eta_{j+1}+M}^{\eta_{j-1}-M}e_j^*(z)\lp f'(\varphi_j(z))-f'(0)\rp \lp\varphi_{j+1}(z) \rp\rmd z+\textrm{hot} \nonumber\\ 
I_2    = & I_2^\textrm{R}+I_2^\textrm{L}+\textrm{hot}.\nonumber
\end{align}
From \eqref{eq:A5} and the definition of $e^*_j,$ we infer that 
\begin{align}
I_2^\textrm{L} = & e^{\lambda(\eta_j-\eta_{j+1})}\Lambda a_+\int_\R e^*(z)\lp f'(\varphi (z)-f'(0))\rp e^{\lambda z}\rmd z+\textrm{hot} \nonumber \\
= & -a_{\textrm{L}}e^{\lambda(\eta_j-\eta_{j+1})}+\textrm{hot},\label{eq:A7}\\
I_2^\textrm{R} = & -e^{\mu(\eta_j-\eta_{j-1})}\Lambda a_-\int_\R e^*(z)\lp f'(\varphi(z))-f'(0)\rp e^{\mu z}\rmd z+\textrm{hot} \nonumber \\
=& a_{\textrm{R}}e^{\mu(\eta_j-\eta_{j-1})}+\textrm{hot}, \label{eq:A8}
\end{align}
where the coefficients $a_{\textrm{L}},a_{\textrm{R}}$ are precisely those given in \eqref{eq:20}. To complete 
our derivation of \eqref{eq:20} it remains to determine the sign of these coefficients. 
We present the computation for $a_{\textrm{L}},$ the case $a_{\textrm{R}}$ being similar. It uses estimates 
\eqref{eq:41} for $\varphi$ and its derivatives up to order 2, as well as the observation 
that $\ds \varphi'f'(\varphi)=-\varphi'''-c\varphi''.$ It follows that 
\begin{align}
a_{\textrm{L}}= & -\Lambda a_+\int_\R e^{cz}\varphi'(z)\lp f'(\varphi (z)-f'(0))\rp e^{\lambda z}\rmd z \nonumber \\
 = & -\Lambda a_+\int_\R e^{-\mu z}\lp -\varphi'''(z)-c\varphi''(z)-f'(0)\varphi'(z)\rp\rmd z \label{eq:A9}\\ 
 = & \Lambda a_+\lp 2\mu^2a_-+c\mu a_-+\lp \mu^2+c\mu+f'(0)\rp\int_\R e^{-\mu z}\varphi'(z)\rmd z\rp \label{eq:A10}
\end{align}
with three successive integration by parts from \eqref{eq:A9} to \eqref{eq:A10}. 
Finally, notice that $\mu^2+c\mu+f'(0)=0,$ and one gets 
\begin{equation}\label{eq:A11}
a_{\textrm{L}}=\Lambda a_+ a_-\mu\lp 2\mu+c\rp>0. 
\end{equation}
A quite similar computation gives
\begin{equation}\label{eq:A12}
a_{\textrm{R}}=\Lambda a_+ a_-\lambda\lp 2\lambda +c\rp>0  
\end{equation}

\paragraph{Negligible terms: integration away from the front} It remains to prove that 
the integral $I_2$ is actually the one that drives the dynamics, and to discuss which value of $M$ would be appropriate.We focus on $I_1^-$ for the computations, the situation being similar for $I_1^+.$ Le us first notice that linearizing at any front, we obtain that
 for all $z\in\R,$ 
 \begin{equation}\label{eq:A13}
f\lp\sum_{i=1}^n\varphi_i(z)\rp-\sum_{i=1}^n f(\varphi_i(z))=O\lp e^{-\frac{\mu}{2}\delta_{\min}}\rp,
 \end{equation}
 while, linearizing around at $\varphi_{j+1},$ we have that for all $z\in[\eta_{j+1},\eta_{j+1}+M],$
 \begin{equation}\label{eq:A132}
f\lp\sum_{i=1}^n\varphi_i(z)\rp-\sum_{i=1}^n f(\varphi_i(z))=O\lp e^{\mu(z-\eta_j)}+e^{\lambda(z-\eta_{j-2})}\rp.
 \end{equation}
 Then, once again we split the integration domain into two.
 \begin{align}
\frac{1}{\Lambda}I_1^- & = \int_{-\infty}^{\eta_{j+1}+M}e_j^*(z) \lp f\lp\sum_{i=1}^n\varphi_i(z)\rp-\sum_{i=1}^n f(\varphi_i(z))\rp \rmd z \nonumber \\
 & = \int_{-\infty}^{\eta_{j+1}} e_j^*\lp f\lp\sum_{i=1}^n\varphi_i\rp-\sum_{i=1}^n f(\varphi_i)\rp \rmd z+ \int_{\eta_{j+1}}^{\eta_{j+1}+M} e_j^*\lp f\lp\sum_{i=1}^n\varphi_i\rp-\sum_{i=1}^n f(\varphi_i)\rp \rmd z \nonumber \\
  & = J_1+J_2. \nonumber
\end{align}
Combining \eqref{eq:A13}, \eqref{eq:41} and the definition of $e^*_j,$ we get that
\begin{align}
J_1 & = \int_{-\infty}^{\eta_{j+1}}e^{c(z-\eta_j)}\varphi'(z-\eta_j)f\lp\sum_{i=1}^n\varphi_i(z)\rp-\sum_{i=1}^n f(\varphi_i(z))\rmd z \nonumber \\
\lb J_1\rb & \leq C\int_{-\infty}^{\eta_{j+1}}e^{c(z-\eta_j)}e^{\mu(z-\eta_j)}\rmd z\  e^{-\frac{\mu}{2}\delta_{\min}}\leq C e^{(\lambda-\mu/2)\delta_{\min}} \label{eq:A14}
\end{align}
for some positive constant $C.$ On the other hand, using \eqref{eq:A132} and \eqref{eq:41},
one has that for some positve constant $C,$ possibly distinct from the above one,
\begin{align}
\lb J_2\rb & \leq C\int_{\eta_{j+1}}^{\eta_{j+1}+M}e^{-\lambda(z-\eta_j)}\lp e^{\mu(z-\eta_j)}+e^{\lambda(z-\eta_{j+2})}\rp\rmd z \nonumber \\
 & \leq C\lp e^{(-2\mu-c)(\delta_j-M)}+Me^{\lambda(\delta_j+\delta_{j-1})}\rp \leq C\lp  e^{(-2\mu-c)(\delta_{\min}-M)}+Me^{2\lambda\delta_{\min}}\rp. \label{eq:A15}
\end{align}
Combining \eqref{eq:A14} and \eqref{eq:A15}, we get, for some positive constant $C,$
\begin{equation}\label{eq:A16}
\lb I_1^-\rb \leq C\lp e^{(\lambda-\mu/2)\delta_{\min}}+e^{(-2\mu-c)(\delta_{\min}-M)}+Me^{2\lambda\delta_{\min}}\rp.
\end{equation}
Similar computations give the following estimates, for some positive constant $C,$
\begin{equation}\label{eq:A17}
\lb I_1^+\rb\leq C\lp e^{-\frac{3}{2}\mu\delta_{\min}}+e^{(\lambda-\mu)(\delta_{\min}-M)}+Me^{-2\mu\delta_{\min}}\rp.
\end{equation}
It remains to choose $M$ and determine $\alpha>0$ such that both $I_1^-,I_1^+,$ and the higher order terms in \eqref{eq:A6} are of order $e^{-(\mu+\alpha)\delta_{\min}},$ as long as 
$\delta_{\min}$ is large enough. Let us fix $\alpha>0$ small enough, and $\varepsilon=\frac{\alpha}{\mu}:$ we fix $\ds M=\lp\frac{1}{2}+\varepsilon\rp\delta_{\min},$
and the estimates are as required. 

\section{Proof of Thm.~\ref{thm:localstability}\label{appendixB} (local stability of $E$)}

We write $\Sigma_n=\sum_{k=1}^n\Psi_k^3-\sum_{k=1}^{n-1}\Psi_k\Psi_{k+1}^2$ so \eqref{angsys2} becomes
\begin{align*}
     \Psi_j'=\begin{cases}\Psi_1\left(\Psi_1^3-\Psi_1-\Psi_1\Psi_2^2+\sum_{k=2}^n\Psi_k^3-\sum_{k=2}^{n-1}\Psi_k\Psi_{k+1}^2\right), & j=1\\
    \Psi_j\left(\Psi_j^3-\Psi_j+\Psi_{j-1}-\Psi_j\Psi_{j+1}^2-\Psi_{j-1}\Psi_j^2+\sum_{k=1,\atop k\neq j}^n\Psi_k^3-\sum_{k=1,\atop k\notin\{j-1,j\}}^{n-1}\Psi_k\Psi_{k+1}^2\right), & 2\leqslant j\leqslant n-1\\
    \Psi_n\left(\Psi_n^3-\Psi_n+\Psi_{n-1}-\Psi_{n-1}\Psi_n^2+\sum_{k=1}^{n-1}\Psi_k^3-\sum_{k=1}^{n-2}\Psi_k\Psi_{k+1}^2\right), & j=n\end{cases}
\end{align*}

In the following, $\chi_A$ denotes the characteristic function which is either $1$ if condition $A$ is satisfied or $0$ otherwise. For instance, for given $i\in\mathbb{N}$, $\chi_{\{i=2\}}=1$ if $i=2$ and $0$ otherwise. The Jacobian $J:=\left(\frac{\partial f_j}{\partial\Psi_i}\right)_{1\leqslant i,j\leqslant n}$ is then given by 
\begin{align*}
    \frac{\partial f_1}{\partial\Psi_i}=\begin{cases}-2\Psi_1+4\Psi_1^3-2\Psi_1\Psi_2^2+\sum_{k=2}^n \Psi_k^3-\sum_{k=2}^{n-1} \Psi_k\Psi_{k+1}^2, & i=1\\
    3\Psi_1\Psi_i^2-2 \chi_{\{i=2\}}\Psi_1^2\Psi_2-\Psi_1\Psi_{i+1}^2-2\chi_{\{i>2\}}\Psi_1\Psi_{i-1}\Psi_i, & 2\leqslant i\leqslant n-1\\
    3\Psi_1\Psi_n^2-2\Psi_1\Psi_{n-1}\Psi_n, & i=n\end{cases}
\end{align*}
and, for $2\leqslant j\leqslant n-1$ and $1\leqslant i\leqslant n$,
\begin{align*}
    \frac{\partial f_j}{\partial\Psi_i}=\begin{cases}\Psi_{j-1}-2\Psi_j+4\Psi_j^3-2\Psi_j\Psi_{j+1}^2-3\Psi_{j-1}\Psi_j^2+\sum_{k=1,\atop k\neq j}^n \Psi_k^3-\sum_{k=1,\atop k\notin\{j-1,j\}}^{n-1}\Psi_k\Psi_{k+1}^2, & i=j\\
    \Psi_j+3\Psi_{j-1}^2\Psi_j-\Psi_{j}^3-2\chi_{\{j>2\}}\Psi_{j-2}\Psi_{j-1}\Psi_j, & i=j-1\\
    3\Psi_j\Psi_{j+1}^2-2\Psi_j^2\Psi_{j+1}-\Psi_j\Psi_{j+2}^2, & i=j+1<n\\
    3\Psi_{n-1}\Psi_n^2-2\Psi_{n-1}^2\Psi_n, & i=j+1=n\\
    3\Psi_j\Psi_i^2-\Psi_j\Psi_{i+1}^2-2\chi_{\{i>1\}}\Psi_j\Psi_{i-1}\Psi_i, & 1\leqslant i\leqslant j-2\\
    3\Psi_j\Psi_i^2-2\Psi_j\Psi_{i-1}\Psi_i-\chi_{\{i<n\}}\Psi_j\Psi_{i+1}^2, & j+2\leqslant i\leqslant n
   \end{cases}
\end{align*}
Finally, for $j=n$ and $1\leqslant i\leqslant n$,
\begin{equation*}
    \frac{\partial f_n}{\partial\Psi_i}=\begin{cases}4\Psi_n^3-2\Psi_n+\Psi_{n-1}-3\Psi_{n-1}\Psi_n^2+\sum_{k=1}^{n-1}\Psi_k^3-\sum_{k=1}^{n-2}\Psi_k\Psi_{k+1}^2, & i=n\\
    \Psi_n-\Psi_n^3+3\Psi_n\Psi_{n-1}^2-2\Psi_{n-2}\Psi_{n-1}\Psi_n, & i=n-1\\
    3\Psi_n\Psi_i^2-\Psi_n\Psi_{i+1}^2-2\chi_{\{i>1\}}\Psi_n\Psi_{i-1}\Psi_i, & \textrm{else}
\end{cases}
\end{equation*}

It remains to determine the eigenvalues of $J(E)$. To this end, since $n$ is fixed throughout, we set $b:=n(n+1)(2n+1)=2n^3+3n^2+n$ and $\beta=\frac{b}{6}$, $\alpha:=\sqrt{1/\beta}$. Then $\alpha=\sqrt{\frac{6}{n(n+1)(2n+1)}}$ so that  $E=\alpha(1,2,\ldots,n)^\intercal$. After some straightforward computations, the first and last row of $J(E)$ are given by 
\begin{align*}
        \frac{\partial f_1}{\partial\Psi_i}(E)=
        \alpha^3\cdot\begin{cases}-\beta-1, & i=1\\
        -1, & 2\leqslant i\leqslant n-1\\
        n^2+2n, & i=n
         \end{cases}
&   & \quad     \frac{\partial f_n}{\partial\Psi_i}(E)=
        \alpha^3\cdot\begin{cases}n(n^2+2n-\beta), & i=n\\
        n(\beta-1), & i=n-1\\
        -n, & \textrm{else}
        \end{cases}
    \end{align*}
    For $2\leqslant j\leqslant n-1$ and $1\leqslant i\leqslant n$, the entries are
    \begin{align*}
        \frac{\partial f_j}{\partial\Psi_i}(E)=\alpha^3\cdot\begin{cases}-j(1+\beta(n)), & i=j\\
       j(\beta(n)-1), & i=j-1\\
        -j, & i=j+1<n\\
        n^3+n^2-2n, & i=j+1=n\\
        -j, & 1\leqslant i\leqslant j-2\\
        jn(n+2), & j+2\leqslant i=n\\
        -j, & j+2\leqslant i<n
        \end{cases}
    \end{align*}

Let $C:=\{\vec{c}_1,\vec{c}_2,\ldots,\vec{c}_n\}:=\alpha^{-3}J(E)$, where $\vec{c}_i$ are the column vectors. First note that the radial direction from $E$, $\vec{v}_1:=(1,2,3,\ldots,n)^\intercal$, is an eigenvector of $C$ with eigenvalue $\lambda_1:=2\beta>0$. However, this unstable direction is transverse to the invariant sphere and thus irrelevant here. Also by invariance, all remaining eigenvectors are orthogonal to $\vec{v}_1$.

A basis of the orthogonal complement $(\textrm{span}(\vec{v}_1))^\perp$ is given by
\begin{align*}
  V&:=\{\vec{v}_2,\vec{v}_3,\ldots,\vec{v}_n\},\quad 
  \vec{v}_j:=(j,0,\ldots,0,-1,0,\ldots,0)^\intercal,\quad j=2,\ldots,n-1,
\end{align*}
where the entry $-1$ in $\vec{v}_j$ is on the $j$-th position. In order to reduce the matrix $C$ onto the tangential dynamics, we orthogonalize the vectors $\vec{v}_j$, $j=1,2,\ldots,n$ by Gram-Schmidt and normalize afterwards. This results in the column vectors of the matrix $W$ given by
\begin{align*}
    \vec{w}_1&:=q_1\vec{v}_1, \quad
    \vec{w_k}:=q_k(1,2,\ldots,k-1,-a_k,0,\ldots,0)^\intercal,\quad k=2,3,\ldots,n,
\end{align*}
with $a_k:=\frac{1}{k}\sum_{i=1}^{k-1}i^2=\frac{1}{6}(k-1)(2k-1)$ and normalizing constants 
\begin{align*}
    q_1&:=\left(\sum_{i=1}^{n}i^2\right)^{-1/2}=\frac{\sqrt{6}}{\sqrt{n(n+1)(2n+1)}}, \quad
    q_k:=\left(\sum_{i=1}i^2+a_k^2\right)^{-1/2}=\frac{6}{\sqrt{4k^4-5k^2+1}}.
\end{align*}
We claim that changing the basis with $W$ yields the lower triangular matrix
\begin{equation*}
    W^\intercal CW=\begin{pmatrix}\lambda_1 & 0 & 0 & \ldots & 0\\0 & b_{2,2} & 0 & \ldots & 0\\
    0 & b_{3,2} & b_{3,3} & \ldots & 0\\\vdots & \vdots & \vdots & \ddots & \vdots\\
  0 & b_{n,2} & b_{n,3} & \ldots & b_{n,n}.\end{pmatrix}
\end{equation*}
with diagonal entries 
\begin{equation*}
    b_{k,k}=-k\beta=-\frac{k}{6}n(n+1)(2n+1)<0.
\end{equation*}

In order to prove this, we first note that with 
\begin{equation*}
\vec{d}_k:=(\vec{w}_k\cdot\vec{c}_1,\vec{w}_k\cdot \vec{c}_2,\ldots,\vec{w}_k\cdot\vec{c}_n)^\intercal
\end{equation*}
we have $b_{k,k}=\vec{d}_k\cdot\vec{w}_k$ so it remains to determine $\vec{d}_k$ for each $k\in\{2,3,\ldots,n\}$. 

For $2\leqslant k\leqslant n-1$ and $1\leqslant j\leqslant n-1$,
\begin{align*}
 \vec{w}_k\cdot \vec{c}_j=\begin{cases}-q_k\sum_{i=1}^{k-1}i^2+ka_kq_k, & k<j\\
 -q_k\sum_{i=1}^{k-1}i^2+ka_kq_k(1+\beta), & k=j\\
 -q_k\sum_{i=1}^{j-1}i^2-j^2q_k(1+\beta)-ka_kq_k(\beta-1), & k=j+1\\
 -q_k\sum_{\substack{1\leqslant i\leqslant k-1,\\ i\notin\{j,j+1\}}}i^2-j^2q_k(1+\beta)+ka_kq_k+(j+1)^2q_k(\beta-1), & k\geqslant j+2
 \end{cases}   
 \end{align*}
For $2\leqslant k\leqslant n-1$, we have $\vec{w}_k\cdot\vec{c}_n=q_kn(n+2)\left(\sum_{i=1}^{k-1}i^2-ka_k\right)$; finally, for $j<n$,
\begin{equation*}
    \vec{w}_n\cdot\vec{c}_j=-q_n\left(\sum_{\substack{1\leqslant i\leqslant n-1,\\ i\notin\{j,j+1\}}}i^2+j(1+\beta)-(j+1)^2(\beta-1)+a_nn(n^2+2n-\beta)\right)
\end{equation*}
and $\vec{w}_n\cdot\vec{c}_n=q_n\left(n(n+2)\sum_{i=1}^{n-1}i^2-a_nn(n^2+2n-\beta)\right)$. 

This leads to the form $b_{k,k}= A_1 + A_2 + A_3$ with (after straightforward computations)
\begin{align*}
  &A_1:=\sum_{j=1}^{k-2}\left[jq_k^2\left(ka_k+(j+1)^2(\beta-1)-j^2(1+\beta)-\sum_{\substack{1\leqslant i\leqslant k-1,\\ i\notin\{j,j+1\}}}i^2\right)\right]=\left(\frac{3}{2k-1}+\frac{42}{k+1}-\frac{75}{2k+1}\right) \beta,\\
  &A_2:=-(k-1)q_k^2\left(\sum_{i=1}^{k-2}i^2+(k-1)^2(1+\beta)+ka_k(\beta-1)\right)=\left(\frac{72}{2k+1}-\frac{36}{k+1}-\frac{3}{2k-1}-3\right)\beta,\\
  &A_3:=a_kq_k^2\sum_{i=1}^{k-1}i^2-ka_k^2q_k^2(1+\beta)=\left(\frac{3}{2k+1}+3-k-\frac{6}{k+1}\right)\beta,
\end{align*}
and thus $b_{k,k}=-k\beta$. We conclude that the eigenvalues  of $C$ are
\begin{align*}
    \lambda_{1,2}&=\pm 2\beta, \quad
    \lambda_{k}=-k\beta,\quad 3\leqslant k\leqslant n.
\end{align*}
and thus $E$ is locally exponentially stable on $S^{n-1}$.

\section{Implementation of the simulations}\label{app:sim}
In this section, we briefly describe the numerical methods and software used in the illustrating figures. The software packages we use are \textsc{JCASim}  \cite{inproceedings} for simulating cellular automata, \textsc{pde2path} \cite{NMTMA-7-58} for PDE simulation, continuation and bifurcation package, and \textsc{Mathematica}.

\medskip 
The space-time plot of the Greenberg-Hastings cellular automaton (Fig.~\ref{fig:GHCATheta}~(a)) has been created with \textsc{JCASim} via an adaption of the method \texttt{transition}, where we specify the transition from one generation of cell states to the next according to the Greenberg-Hastings rule for $e=2$ excited and $r=4$ refractory states.

The phase portraits on the spheres (Fig.~\ref{fig:spheres}) have been created with \textsc{Mathematica} using the standard routines \texttt{StreamPlot} and \texttt{TransformedField} to change from Cartesian to spherical coordinates.

For the simulations of the PDE~\eqref{theta-eq} plotted in \ref{fig:GHCATheta}~(b), Fig.~\ref{fig:ordering}, \ref{fig:washingout} and \ref{fig:kinkdistPBC} , we have used \textsc{pde2path} 
with $\muf=0.05$ or $\muf=0.2$, respectively. As discretization of the interval \texttt{(-lx,lx)} by mesh width \texttt{2*lx/dsc} we have chosen \texttt{lx=10} and \texttt{dsc=500} which corresponds to a mesh width of $h=0.04$. Initial data \texttt{p.u(1:p.nu)} have been step functions with equilibrium levels at \texttt{eq+2*pi*j} for $j=1,2,\ldots,k, k\in\mathbb{N}$. 

For the time stepping the \textsc{pde2path} routines \texttt{tints} or \texttt{tintsfreeze} have been used since most simulations go over a rather long time period in order to catch exponentially small interaction effects between kinks (antikinks). \texttt{tintsfreeze} removes the translational symmetry so that the required spatial interval need not be excessively long. In particular, since annihilations of kinks and antikins on the one side and interactions of kinks (antikinks) on the other side happen at different time scales (the annihilations being much faster), it is nearly impossible to observe both aspects simultaneously (even when choosing very large distances or a large space interval). We have circumvented this problem in Fig.~\ref{fig:washingout} by starting with pure kink data ``sending in`` antikinks from the left space boundary repeatedly after time has allowed for observable exponentially small effects between the kinks that are left after each annihilation. This has been done by calling \texttt{tintsfreeze} multiple times and using the perturbed result as initial condition for the next call.

The detection of the kink (antikink) positions has been done via a function that stores the (geometrically defined) positions for each time step. 
The differences between these positions are the basis for the distance plots in Fig.~\ref{fig:ordering}, \ref{fig:washingout} and \ref{fig:kinkdistPBC}. 
For better illustration we have smoothened the distance plots by averaging the data and interpolating by splines using the \texttt{MATLAB} routine \texttt{interp1}. 
Finally, for general information about \textsc{pde2path} and, in particular, the simulation with periodic boundary conditions (Fig.~\ref{fig:kinkdistPBC}), we refer to \cite{enoc2014} and the \textsc{pde2path} webpage, respectively.

\def\cprime{$'$}
\bibliographystyle{alpha}


\end{document}